\documentclass{article}

\usepackage[ruled,vlined,linesnumbered,commentsnumbered]{algorithm2e}

\usepackage{graphicx}

\usepackage{authblk}

\usepackage{tikz}
\usetikzlibrary{arrows,calc}

\tikzstyle{arc}=[->,shorten <=3pt, shorten >=3pt, >=stealth, line width=1.1pt]
\tikzstyle{edge}=[shorten <=2pt, shorten >=2pt, >=stealth, line width=1.1pt]
\tikzstyle{gryE}=[shorten <=2pt, shorten >=2pt, >=stealth, line width=1.1pt,
color=gray!60]
\tikzstyle{vertex}=[circle, fill=white, draw, minimum size=5pt, inner sep=0pt,
outer sep=0pt,line width=0.7pt]
\tikzstyle{blackV}=[circle, fill=black, draw, minimum size=5pt, inner sep=0pt,
outer sep=0pt]
\tikzstyle{grayV}=[circle, fill=white, draw=gray!65, minimum size=5pt, inner
sep=0pt, outer sep=0pt, line width=1]

\usepackage{float}

\usepackage{xcolor}

\usepackage{amsmath}
\usepackage{amssymb}
\usepackage{stmaryrd}

\usepackage{complexity}

\usepackage{amsthm}

\usepackage{thmtools}

\usepackage{enumerate}

\usepackage[colorlinks=true]{hyperref}

\usepackage[capitalize]{cleveref}

\usepackage{xspace}

\newtheorem{theorem}{Theorem}[section]
\newtheorem{lemma}[theorem]{Lemma}
\newtheorem{corollary}[theorem]{Corollary}
\newtheorem{proposition}[theorem]{Proposition}
\newtheorem{problem}[theorem]{Problem}
\newtheorem{redrule}{Rule}[section]

\usepackage{complexity}

\usepackage{multirow}
\usepackage{tabularx}
\usepackage{fbox}
\usepackage{hhline}
\newcommand{\defproblem}[3]{\par
 \vspace{1mm}
\noindent
\begin{center}
\fbox{
 \begin{minipage}{0.7\textwidth}
 \begin{tabular*}{\textwidth}{@{\extracolsep{\fill}}lr} #1 &  \vspace{1mm} \\ \end{tabular*}
{\textbf{Input:}} #2
  \vspace{1mm}\\%
 {\textbf{Question:}} #3
 \end{minipage}
 }
 \end{center}
 \vspace{1mm}\par
}

\def\mayatupi{Maya-Tupi\xspace}
\def\mtrecognition{{\sc \MT-recognition}\xspace}
\def\MT{MT\xspace}
\def\barM{\overline{M}\xspace}
\def\no{\texttt{no}}
\def\yes{\texttt{yes}}

\newcommand{\nd}{{\sf nd}\xspace}

\newcommand{\tw}{{\sf tw}\xspace}
\newcommand{\cw}{{\sf cw}\xspace}

\newcommand{\interval}[1]{\llbracket #1 \rrbracket}

\DeclareMathOperator{\Fdisc}{\ensuremath{\mathcal{F}_{disc}}}
\DeclareMathOperator{\Fcog}{\ensuremath{\mathcal{F}_{cog}}}

\title{Maya-Tupi graphs: a generalization of split graphs%
}

\author[1]{J\'ulio~Ara\'ujo\thanks{julio@mat.ufc.br}}
\author[2]{C\'esar~Hern\'andez-Cruz\thanks{chc@ciencias.unam.mx}}
\author[3]{Cl\'audia~Linhares~Sales\thanks{linhares@dc.ufc.br}}

\affil[1]{Departamento de Matem\'atica,  Universidade Federal do Cear\'a, Fortaleza, Brazil}
\affil[2]{Departamento de Matem\'aticas, Facultad de Ciencias, Universidad Nacional Autónoma de M\'exico, M\'exico}
\affil[3]{Departamento de Computa\c c\~ao,  Universidade Federal do Cear\'a, Fortaleza, Brazil}

\begin{document}
\date{}

\maketitle
\begin{abstract}
  We define the family of \mayatupi graphs as those graphs that admit a
  partition $(A,B)$ of their vertex sets such that $A$ induces a complete
  multipartite graph where each part has size at most two, and $B$ induces a
  graph where every connected component is $K_1$ or $K_2$.   The family of
  \mayatupi graphs is self complementary, generalizes split graphs, falls into
  the sparse-dense partitioning schema and is characterized by finitely many
  forbidden induced subgraphs.   Unfortunately, our computational experiments
  show that the number of minimal forbidden induced subgraphs to characterize
  \mayatupi graphs is greater than 2000.

  In this work, we find characterizations in terms of minimal forbidden induced
  subgraphs for disconnected graphs, which imply the same for cographs; our
  results imply linear-time certifying recognition algorithms for \mayatupi
  graphs within these classes. We also show that \mayatupi graphs can be
  recognized in $\mathcal{O}(n^3)$-time in $C_4$-free graphs and in graphs with
  bounded neighborhood diversity; in $\mathcal{O}(n^4)$-time for triangle-free
  graphs; and in $\mathcal{O}(n^2)$-time for graphs with bounded clique-width.

  We provide efficient algorithms to calculate the clique, the independence, the
  chromatic, and the treewidth numbers, as well as a minimum fill-in for
  \mayatupi graphs.
\end{abstract}


\section{Introduction}

For basic definitions on Graph Theory and Computational Complexity, we refer
to~\cite{bondy2008,sipser1996}. All graphs in this work are simple, finite, with
neither loops nor parallel edges.  When we consider a \textit{partition} of the
vertex set of a graph, we allow parts to be empty.

A graph $G$ is \emph{split} if its vertex set can be partitioned into an
independent set $S$ and a clique $K$. There is a vast literature on split
graphs, not only presenting properties about their structure, but also
concerning the study of graph problems restricted to this
class~\cite{merris2003}.   Among the many generalizations of split graphs, some
of the better known are $(k,l)$-graphs, defined by Brandst\"adt
\cite{brandstadtDM152} as graphs that admit a partition $(A,B)$ of its vertex
set such that $A$ induces a $k$-colorable graph and $B$ induces the complement
of an $l$-colorable graph; the $k$-split graphs, defined by Chudnovski and
Seymour~\cite{CHUDNOVSKY201411} as those admitting a partition $(A,B)$ of their
vertex set such that $A$ induces a graph with clique number at most $k$ and $B$
induces a graph with stability number at most $k$; and $(s,k)$-polar graphs,
introduced by Ekim, Mahadev and de Werra \cite{ekimDAM156} as graphs whose
vertex set admits a partition $(A,B)$ such that $A$ induces a complete
$s$-partite graph and $B$ induces the complement of a complete $k$-partite
graph.   Interestingly, more than two decades before, Tyshkevich and Chernyak
\cite{tyshkevichK2} introduced $(s,k)$-polar graphs as those graphs admitting a
partition of their vertex set into two sets $A$ and $B$, such that the connected
components of $\overline{G[A]}$ and $G[B]$ are complete graphs of size at most
$s$, and $k$, respectively.   A possible explanation for the existence of two
different definitions of $(s,k)$-polar graphs comes from the definition of a
\emph{polar graph}, also introduced in \cite{tyshkevichK2}, as a graph admitting
a vertex partition $(A,B)$ such that $A$ induces a complete multipartite graph
and $B$ the complement of a complete multipartite graph.   Clearly, the
definition of polar graph coincides with the intuitive notion of $(\infty,
\infty)$-polar graphs for both definitions of $(s,k)$-polar graphs.  It seems
that the class of polar graphs lived long enough for the community to forget the
definition of $(s,k)$-polar graphs and introduce a new one.   Henceforth, in
order to propose a solution that permits both definitions to coexist, we use the
term $\langle s,k \rangle$-polar graph to refer to Tyshkevich and Chernyak's
definition.

The partitions for all the families described in the previous paragraph are
\textit{sparse-dense} partitions in the sense of \cite{federSIDMA16}.   Since,
bipartite graphs, complete multipartite graphs, and their complements, can be
recognized in polynomial time, Theorem 3.1 in \cite{federSIDMA16} implies that
$(k,l)$-graphs with $\max \{ k, l \} \le 2$, $(s,k)$-polar graphs, and $\langle
s,k \rangle$-polar graphs are recognizable in polynomial time. In contrast,
Chernyak~\cite{chernyakDM62} proved that the problem of recognizing polar graphs
is \NP-complete, and it is well-known that recognizing $(k,l)$-graphs with $\max
\{ k, l \} \ge 3$ is an \NP-complete problem.   Ekim, Hell and Stacho proved in
\cite{ekimDAM156a} that, when restricted to chordal graphs, it is possible to
recognize polar graphs in polynomial time, even when the number of forbidden
induced chordal subgraphs is infinite for this family.   For cographs, Ekim,
Mahadev and de Werra exhibit in \cite{ekimDAM156} the complete list of forbidden
induced cograph obstructions for polarity, presenting also a recognition
algorithm running in linear time.

In this work, we focus on a generalization of split graphs that we think is of
particular interest, and explore some of its basic properties.   A graph $G$ is
\emph{\mayatupi} if its vertex set can be partitioned into two sets $A$ and $B$
such that $\delta(G[A])\geq |A|-2$ and $\Delta(G[B])\leq 1$. We say that the
pair $(A,B)$ is an \emph{$\MT$-partition} of $G$. Consequently, $G[A]$ is
$\{\overline{P_3},\overline{K_3}\}$-free, while $G[B]$ is $\{P_3,K_3\}$-free. We
can also obtain an equivalent definition related to defective colorings
\cite{cowenJGT10}, as $B$ induces a $(1,1)$-colorable graph, and $A$ induces the
complement of a $(1,1)$-colorable graph.   Additionally, note that \mayatupi
graphs are naturally the $\langle 2,2 \rangle$-polar graphs.   It follows from
the different equivalent definitions that, just like split graphs, the class of
\mayatupi graphs is self-complementary and hereditary.

We introduce a final equivalent definition which will be useful in some proofs.
A graph $G$ is \mayatupi if and only if its vertex set can be partitioned into
four (possibly empty) sets $S$, $M$, $\barM$ and $K$ such that: $S$ is
independent; $M$ induces a matching; $\barM$ induces a complete graph minus a
perfect matching (or an anti-matching); $K$ is a clique; $S$ is anticomplete to
$M$; and $K$ is complete to $\overline{M}$ (notice that when all parts are
non-empty, then this is also a skew partition of $G$). If $(A,B)$ is an
$\MT$-partition of $G$, note that $\{K,\barM\}$ is a partition of $A$, while
$\{S,M\}$ is a partition of $B$. Abusing nomenclature, we also say that
$(K,\barM,S,M)$ is an \emph{$\MT$-partition} of $G$. The number of elements of
the tuple defines to which partition we refer.

\begin{figure}[ht!]
  \centering
  \begin{tikzpicture}
    \begin{scope}[yscale=1,xscale=-1]
    \draw [rounded corners] (-0.3,-2.5) rectangle (1.3,0.7); \node (b) at (-0.6,
    0.5){$B$}; \draw [rounded corners] (2.1,-2.5) rectangle (3.8,0.7); \node (a)
    at (4.1,0.5){$A$};
    
    \node [vertex] (s1) at (0,0){}; \node [vertex] (s2) at (1,0){}; \node
    [vertex] (s3) at (0.5,-1){}; \node [vertex] (s4) at (0.0,-2){}; \node
    [vertex] (s5) at (1.0,-2){};
    
    \draw [edge] (s1) to (s2); \draw [edge] (s4) to (s5);

    \node [blackV] (k1) at (2.5,-1){}; \node [blackV] (k2) at (3.5,0){}; \node
    [blackV] (k3) at (3,-2){}; \draw [edge] (k1) to (k3); \draw [edge] (k2) to
    (k3);

    \draw [gryE] (s1) to [bend left] (k2); \draw [gryE] (s2) to (k1); \draw
    [gryE] (s3) to (k1); \draw [gryE] (s3) to (k2); \draw [gryE] (s5) to (k3);
  \end{scope}
  \begin{scope}[xscale=-1, xshift=-6cm]
    \draw [rounded corners] (-0.3,-2.5) -- (1.3,-2.5) -- (1.3,-1.5) --
        (1.3,0.7) -- (-0.3, 0.7) -- (-0.3,-0.3) -- (1.1,-0.3) -- (1.1,-1.7)
        -- (-0.3,-1.7) -- cycle;
    \node (b) at (-0.6, 0.5){$M$};
    \draw [rounded corners] (-0.2,-1.5) rectangle (0.9,-0.5);
    \node (b) at (-0.5, -0.7){$S$};
    \draw [rounded corners] (2.1,-1.3) rectangle (3.8,0.7);
    \node (a) at (4.1,0.5){$\overline{M}$};
    \draw [rounded corners] (2.1,-1.5) rectangle (3.8,-2.5);
    \node (a) at (4.1,-1.65){$K$};
    
    \node [vertex] (s1) at (0,0){};
    \node [vertex] (s2) at (1,0){};
    \node [vertex] (s3) at (0.5,-1){};
    \node [vertex] (s4) at (0.0,-2){};
    \node [vertex] (s5) at (1.0,-2){};
    
    \draw [edge] (s1) to (s2);
    \draw [edge] (s4) to (s5);

    \node [blackV] (k1) at (2.5,-1){};
    \node [blackV] (k2) at (3.5,0){};
    \node [blackV] (k3) at (3,-2){};
    \draw [edge] (k1) to (k3);
    \draw [edge] (k2) to (k3);

    \draw [gryE] (s1) to [bend left] (k2);
    \draw [gryE] (s2) to (k1);
    \draw [gryE] (s3) to (k1);
    \draw [gryE] (s3) to (k2);
    \draw [gryE] (s5) to (k3);
  \end{scope}
  \end{tikzpicture}
  \caption{An example of a \mayatupi graph containing induced copies of $C_4$,
  $2K_2$ and $C_5$.}
\label{fig:firstexample}
\end{figure}
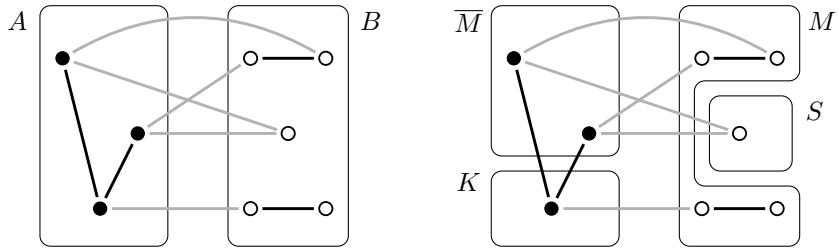

An example of an $\MT$-graph is shown in \cref{fig:firstexample} with $A$
represented by black vertices, $B$ represented by white vertices, and edges
between $A$ and $B$ in gray.   On the right we present the same graph but with
the partition corresponding to the equivalent definition introduced in the
previous paragraph.   It is well-known that $G$ is split if and only if it is
$\{C_4,C_5,2K_2\}$-free~\cite{foldesCNXIX}; note that these three induced
subgraphs occur in the graph depicted in \cref{fig:firstexample}.

Gagarin and Metelsky~\cite{gagarinVNANBSFMN3} proved that there are finitely
many minimal $\langle 1,2 \rangle$-polar obstructions and exhibited the complete
list of $18$ such obstructions. Later, Zverovich and Zverovich
\cite{zverovichGTNNY44} proved that for any choice of $s$ and $k$, there are
finitely many minimal obstructions for the class of $\langle s,k \rangle$-polar
graphs.   As noted above, $\langle 1,1 \rangle$-tc graphs (split graphs)
coincide with $\{ 2K_2, C_4, C_5 \}$-free graphs.  We ran a computer search
which showed that \mayatupi graphs have over 2000 minimal forbidden induced
subgraphs. Therefore, \mayatupi graphs represent a very interesting challenge in
terms of characterization by forbidden induced subgraphs.

The rest of this work is organized as follows. In \Cref{sec:prelim} we introduce
the basic notation that will be used in the rest of the article. The problem of
recognizing \mayatupi graphs is addressed in \Cref{sec:recognition}, where we
verify that a $\mathcal{O}(n^8)$ solution exists, and show that this time can be
improved when the input graph has bounded clique or independence number, or when
$G$ has bounded cliquewidth, or when $G$ has bounded neighborhood diversity.
\mayatupi cographs are characterized in terms of forbidden induced subgraphs in
\Cref{sec:mt-cographs}; as a necessary intermediate step, in this section we
also characterize the disconnected forbidden induced subgraphs for \mayatupi
graphs.      \Cref{sec:c4-free} is devoted to present an algorithm for
recognizing $C_4$-free \mayatupi graphs in time $\mathcal{O}(n^3)$; as a
necessary ingredient for our algorithm we also provide a certifying linear-time
algorithm to recognize $\langle 1,2 \rangle$-polar graphs.   In
\Cref{sec:optimization} we show that a maximum clique, a maximum independent
set, an optimal coloring, an optimal clique cover, an optimal tree decomposition
and a minimum fill-in can be obtained in polynomial time. Conclusions and open
problems are presented in \Cref{sec:conclusions}.


\section{Preliminaries}
\label{sec:prelim}

All graphs are finite, with neither loops nor parallel arcs.   Unless otherwise
stated, when we refer to a graph $G$ we will assume that it has vertex set $V$
and edge set $E$.   Hence, we will often describe the complexity of an algorithm
in terms of $|V|$ and $|E|$.

For a fixed graph $H$, a graph $G$ is \emph{$H$-free} if it has no induced
subgraph isomorphic to $H$.  For a family of graphs $\mathcal{F}$, $G$ is
\emph{$\mathcal{F}$-free} if it is $H$-free for every $H \in \mathcal{F}$.

Recall that a graph property $\mathcal{P}$ is \emph{hereditary} if it is closed
under taking induced subgraphs, i.e., whenever $G \in \mathcal{P}$ and $H$ is an
induced subgraph of $G$, we also have $H \in \mathcal{P}$.  It follows from the
definition that if $G \notin \mathcal{P}$, then $G$ contains an induced subgraph
$H$ which is also not in $\mathcal{P}$ but such that each of its proper induced
subgraphs is in $\mathcal{P}$; we call $G$ a \textit{$\mathcal{P}$-obstruction}
and $H$ a \textit{minimal $\mathcal{P}$-obstruction}.  In particular, when
$\mathcal{P}$ is the class of \mayatupi graphs, we use the terms
\textit{\MT-obstruction} and \textit{minimal \MT-obstruction}.

We usually assume that graphs have $n$ vertices and $m$ edges. Given graphs $G$
and $H$ with disjoint sets of vertices, the \emph{disjoint union} of $G$ and $H$
is the graph $G+H$ such that $V(G+H) = V(G)\cup V(H)$ and $E(G+H) = E(G)\cup
E(H)$. For a non-negative integer $p$, we denote by $pG$ the graph obtained by
the disjoint union of $p$ copies of $G$.

In graph $G$ with $S$ and $T$ being subsets of $V(G)$, we say that $S$ is
\textit{complete} (\textit{anticomplete}) to $T$ if, for every pair of vertices
$s \in S$ and $t \in T$, $st$ is (resp. not) an edge of $G$.

For a positive integer $k$, we use $\interval{k}$ to denote the set $\{ 1,
\dots, k \}$.


\section{Recognition algorithms}
\label{sec:recognition}

As we mentioned earlier, Zverovich and Zverovich \cite{zverovichGTNNY44} proved
that for any choice of $s$ and $k$, there are finitely many minimal obstructions
for the class of $\langle s,k \rangle$-polar graphs. Since the number of minimal
obstructions to the class of \mayatupi graphs is finite, as they are the
$\langle 2,2 \rangle$-polar graphs, there is a brute-force polynomial-time
algorithm to solve the following decision problem:

\defproblem{\mtrecognition}
{A graph $G$.}
{Is $G$ \mayatupi?}

With computer assistance, we found over 2000 minimal \MT-obstructions having
from 7 to 9 vertices, none with at most 6 or exactly 10 vertices. A natural
question is:

\begin{problem}\label{prob:obstructionwith11}
Does there exist a minimal \MT-obstruction with at least 11 vertices?
\end{problem}

If we answer \cref{prob:obstructionwith11} negatively, this would imply in an
$\mathcal{O}(n^9)$-time algorithm to solve \mtrecognition which does not
construct an \MT-partition if the input graph is indeed a \mayatupi graph, but
returns an \MT-minimal obstruction as a \no-certificate when it is not.   The
fact that this algorithm does not construct an \MT-partition for \yes-instances
of \mtrecognition is unfortunate, as such a partition is the main ingredient of
the algorithms we present in \Cref{sec:optimization}.   Luckily, the work of
Feder, Hell, Klein and Motwani \cite{federSIDMA16} provides an algorithm that
constructs an \MT-partition when the input graph is a \mayatupi graph (although
it does not returns \no-certificates when it is not).

Let $\mathcal{S}$ be the class of graphs having maximum degree $1$, and let
$\mathcal{D}$ be the class of graphs that are complements of graphs in
$\mathcal{S}$.   Clearly, if $G$ is a \mayatupi graph with \MT-partition
$(A,B)$, then $G[A] \in \mathcal{S}$ and $G[B] \in \mathcal{D}$.

\begin{proposition}
\label{pro:mt-recognition}
    There is an algorithm that solves \mtrecognition in time $\mathcal{O}(n^8)$
    and outputs an \MT-partition when one exists.
\end{proposition}
\begin{proof}
    Since both $\mathcal{S}$ and $\mathcal{D}$ are hereditary families and the
    intersection $S \cap D$ has at most $2$ vertices for any $S \in \mathcal{S}$
    and any $D \in \mathcal{D}$, the conclusion is immediate from Theorem 3.1 in
    \cite{federSIDMA16}.
\end{proof}

\Cref{pro:mt-recognition} is important because having a way of constructing the
\MT-partition of a \mayatupi graph enables us to use it to solve a variety of
optimization problems efficiently.   We will do this in \Cref{sec:optimization}.

For the remainder of this section we propose better algorithms to solve
\mtrecognition under the hypothesis that different parameters of the input graph
are bounded by a constant. 

A simple, but powerful, observation is the following. 

\begin{proposition}
  \label{prop:cliqueindependencenumbers}
  Given a graph $G$, if $\alpha(G) \leq k$ or $\omega(G) \leq k$ for some fixed
  positive integer $k$, then \mtrecognition can be solved in time ${\cal
  O}(n^{2k+2})$.
\end{proposition}
\begin{proof}
  Suppose $\omega(G)\leq k$, as the case when $\alpha(G)\leq k$ is analogous.
  One should notice that if $G$ is \mayatupi with \MT-partition $(A,B)$, then
  $|A|\leq 2k$. Thus, we can check all subsets of $V(G)$ of size at most $2k$ to
  be $A$ in an \MT-partition $(A,B)$ of $G$. For each such subset, we can verify
  in linear time whether $B = V(G) \setminus A$ induces a graph with maximum
  degree at most $1$ and whether $A$ induces a complete multipartite graph with
  parts of size at most $2$. Since there are $\mathcal{O}(n^{2k})$ such subsets,
  the result follows.
\end{proof}

\begin{corollary}\label{cor:triaglefreerecognition}
  If $G$ is triangle-free, then \mtrecognition can be solved in time
  $\mathcal{O}(n^4)$.
\end{corollary}

One should also notice that the algorithm described in the proof of
\cref{prop:cliqueindependencenumbers} also constructs an \MT-partition when the
input graph is indeed a \mayatupi graph, but it does not return a
\no-certificate when it is not.

There are many graph parameters that bound the clique or the independence
numbers of a graph such as, for example, chromatic number, degeneracy,
treewidth, vertex cover, feedback vertex set, clique cover number, etc. Thus,
\cref{prop:cliqueindependencenumbers} implies that \mtrecognition can be solved
in polynomial time when the input graph has at least one of these parameters
bounded by a constant. 

A well-known parameter that does not bound neither the clique nor the
independence number is the cliquewidth of a graph $G$, denoted by $\cw(G)$.
However, as we show next, \mtrecognition can also be solved in polynomial time
when the input graph has bounded cliquewidth.

The \emph{cliquewidth} of a graph $G$ is the minimum number of labels needed to
construct $G$ using the following four operations: (1) creation of a new vertex
with a label $i$; (2) disjoint union of two labeled graphs; (3) joining by an
edge every vertex with label $i$ to every vertex with label $j$, for $i \neq j$;
and (4) renaming label $i$ to $j$. A \emph{$k$-expression} of a graph $G$ is an
expression that uses at most $k$ labels and constructs $G$ using the four
operations above.

\begin{proposition}\label{prop:cliquewidth}
  For a fixed positive integer $k$, \mtrecognition can be solved in time
  $\mathcal{O}(n^2)$ if $\cw(G)\leq k$ and the input graph $G$ has $n$ vertices.
\end{proposition}
\begin{proof}
  As proved by Courcelle, Makowsky and Rotics~\cite{Courcelle2000}, any graph
  property definable in monadic second order logic $MSO_1$ can be decided in
  linear time on graphs of bounded cliquewidth, provided that a $k$-expression
  for the input graph is given as part of the input. 
  
  Fomin and Korhonen~\cite{FK2022} presented an algorithm that, for a fixed
  integer $k$, in time $(2^{2^{\mathcal{O}(k)}}) \cdot \mathcal{O}(n^2)$, either
  constructs a $(2^{2k+1}-1)$-expression or correctly concludes that $\cw(G) >
  k$. Thus, if $G$ has cliquewidth at most $k$, for a fixed positive integer
  $k$, we can obtain an $f(k)$-expression of $G$ in time $\mathcal{O}(n^2)$.
  
  It is easy to see that the property of being \mayatupi can be expressed in the
  monadic second order logic $MSO_1$, when it is required that only vertex sets
  are quantified. It suffices to state that there exists a partition of the
  vertex set into two sets $A$ and $B$ such that each vertex in $A$ is adjacent
  to all vertices in $A$ or is not adjacent to exactly one vertex in $A$, and
  that the similar property holds for the vertices in $B$, but with reversed
  adjacency.
  
  Thus, we can decide whether $G$ is \mayatupi in time $\mathcal{O}(n^2)$, if
  $G$ has bounded cliquewidth.
\end{proof}

It is well-known that the algorithm embedded in the proof of Courcelle, Makowsky
and Rotics~\cite{Courcelle2000} is highly exponential in $k$, provided a
$k$-expression of the input graph. To obtain an $f(k)$-expression in quadratic
time, we have two extra layers of exponential in $k$ by using the result of
Fomin and Korhonen~\cite{FK2022}. 

In the sequel, we find a more efficient algorithm to solve \mtrecognition when
the input graph has bounded neighborhood diversity.

Two vertices $u,v$ of a graph $G$ have the same \textit{type} if $N_G(u)
\setminus \{v\} = N_G(v) \setminus \{u\}$. Two vertices with the same type are
\textit{true twins} (resp. \textit{false twins}) if they are adjacent (resp.
non-adjacent). The \textit{neighborhood diversity} of a graph $G$, denoted by
$\nd(G)$, as defined by Lampis~\cite{Lampis12}, is the minimum integer~$w$ such
that $V(G)$ can be partitioned into $w$ sets such that all the vertices in each
set have the same type. Note that the property of having the same type is an
equivalence relation, and that all the vertices in a given type are either true
or false twins, hence defining either a clique or an independent set.

Note that having bounded neighborhood diversity does not bound neither the
clique nor the independence number of a graph. Moreover, observe that a graph
that has bounded neighborhood diversity also has bounded cliquewidth, thus
\cref{prop:cliquewidth} also implies that \mtrecognition can be solved in ${\cal
O}(n^2)$-time when the input graph has bounded neighborhood diversity, but with
a highly exponential cost on the parameter. We improve this running time on $k$
next, but at a greater cost on the input size.

\begin{redrule}\label{rule:fourtwins}
  If a graph $G$ contains a subset $S\subseteq V(G)$ of false (resp. true) twins
  and $|S|\geq 4$, then remove $|S|-3$ vertices of $S$ from $G$.
\end{redrule}

\begin{proposition}\label{prop:fourtwins}
  If $G'$ is obtained from the application of \cref{rule:fourtwins} to $G$, then
  $G$ is \mayatupi if and only if $G'$ is \mayatupi.
\end{proposition}
\begin{proof}
  If $G$ is \mayatupi, then $G'$ is also \mayatupi as any induced subgraph of
  $G$ also belongs to the class.
  
  Suppose that $G'$ is \mayatupi and let $(A,B)$ be an \MT-partition of $G'$.
  Let $S$ be the subset of twins of $G$ with at least four elements such that
  the removal of all but three of these vertices created $G'$. Without loss of
  generality, assume that $S$ is composed of \emph{false} twins as the proof for
  true twins is analogous to what follows. 
  
  Since $S$ is composed of false twins, note that at most two vertices in $S\cap
  V(G')$ lie in $A$. Consequently, at least one vertex in $S\cap V(G')$ lies in
  $B$. Thus, the vertices in $S\setminus V(G')$ can safely be added to $B$ in
  order to obtain an \MT-partition to $G$ as the vertices that will be added
  already have a false twin in $B$ in a valid \MT-partition $(A,B)$ of $G'$. So
  $G$ is also \mayatupi.
\end{proof}

\begin{proposition}\label{prop:nd}
  For a fixed positive integer $k$, \mtrecognition can be solved in time
  $\mathcal{O}\left(n^3+2^{(k+2)\log k}\right)$ if $\nd(G)\leq k$ and the input
  graph $G$ has $n$ vertices.
\end{proposition}
\begin{proof}
  As argued by Lampis~\cite{Lampis12}, the neighborhood diversity of a graph $G$
  can be computed in $\mathcal{O}(n^3)$-time, where $n$ is the number of
  vertices of $G$. In fact, it also provides a partition of $G$ into the
  equivalence classes of types having at most $k$ parts, each of them
  corresponding to a set of true or false twins, by just checking whether all
  pairs of vertices have the same type.

  Provided with such a partition, we can, in linear time, successively apply
  \cref{prop:fourtwins} to remove all but three vertices of each set of false or
  true twins. The resulting graph $G'$ has at most $3k$ vertices, and thus the
  application of a brute-force algorithm to check whether $G'$ is \mayatupi is
  done roughly in ${\cal O}(k^{k+2})$-time. The complete algorithm has thus a
  running time of $\mathcal{O}(n^3 + k^{k+2})$.
\end{proof}

One should notice that the algorithm described in the proof of \cref{prop:nd}
also constructs an \MT-partition when the input graph is indeed a \mayatupi
graph, as the proof of \cref{prop:fourtwins} also explains how to rebuild the
\MT-partition of the original graph from the \MT-partition of the reduced graph,
but it does not return a \no-certificate when it is not.


\section{\mayatupi cographs}
\label{sec:mt-cographs}

\begin{figure}[ht!]
  \centering
  \begin{tikzpicture}
    \begin{scope}[]
      \foreach \i in {0,1,2,3}
      \node [vertex] (\i) at ({90*\i-45}:1){};
      \foreach \i in {0,1,2,3}
      \draw [edge] let \n1 = {int(mod(\i+1,4))}
      in
      (\i) to (\n1);
      \node (label) at (0,-1.3){$C_4$};
    \end{scope}
    
    \begin{scope}[xshift=2.5cm]
      \foreach \i in {0,1,2,3}
      \node [vertex] (\i) at ({90*\i-45}:1){};
      \foreach \i in {0,1,2,3}
      \draw [edge] let \n1 = {int(mod(\i+1,4))}
      in
      (\i) to (\n1);
      \draw [edge] (1) to (3);
      \node (label) at (0,-1.3){diamond};
    \end{scope}
    
    \begin{scope}[xshift=5cm]
      \foreach \i in {0,1,2,3}
      \node [vertex] (\i) at ({90*\i-45}:1){};
      \foreach \i in {0,1,2,3}
      \draw [edge] let \n1 = {int(mod(\i+1,4))}
      in
      (\i) to (\n1);
      \draw [edge] (1) to (3);
      \foreach \i/\j in {0/2,1/3}
      \draw [edge] (\i) to (\j);
      \node (label) at (0,-1.3){$K_4$};
    \end{scope}
    
    \begin{scope}[xshift=7.5cm]
      \foreach \i in {0,...,4}
      \node [vertex] (\i) at ({72*\i+90}:1){};
      \foreach \i in {0,...,4}
      \draw [edge] let \n1 = {int(mod(\i+1,5))}
      in
      (\i) to (\n1);
      \node (label) at (0,-1.3){$C_5$};
    \end{scope}
    
    \begin{scope}[yshift=-3cm]
      \foreach \i in {0,...,4}
      \node [vertex] (\i) at ({72*\i-90}:1){};
      \foreach \i/\j in {0/1,1/2,1/4,3/4,4/0}
      \draw [edge] (\i) to (\j);
      \node (label) at (0,-1.3){bull};
    \end{scope}
    
    \begin{scope}[xshift=2.5cm,yshift=-3cm]
      \foreach \i in {0,1,2,3}
      \node [vertex] (\i) at ({90*\i-45}:1){};
      \node [vertex] (4) at (-0.3,0){};
      \node [vertex] (5) at (0.3,0){};
      \foreach \j in {2,3,5}
      \draw [edge] (4) to (\j);
      \foreach \j in {0,1}
      \draw [edge] (5) to (\j);
      \node (label) at (0,-1.3){$H$};
    \end{scope}
    
    \begin{scope}[xshift=5cm,yshift=-3cm]
      \foreach \i in {0,1,2,3}
      \node [vertex] (\i) at ({90*\i-45}:1){};
      \node [vertex] (4) at (-0.3,0){};
      \node [vertex] (5) at (0.3,0){};
      \foreach \j in {2,3,5}
      \draw [edge] (4) to (\j);
      \foreach \j in {0,1}
      \draw [edge] (5) to (\j);
      \draw [edge] (0) to (1);
      \node (label) at (0,-1.3){$X_{95}$};
    \end{scope}
    
    \begin{scope}[xshift=7.5cm,yshift=-3cm]
      \foreach \i in {0,1,2,3}
      \node [vertex] (\i) at ({90*\i-45}:1){};
      \node [vertex] (4) at (-0.3,0){};
      \node [vertex] (5) at (0.3,0){};
      \foreach \j in {2,3,5}
      \draw [edge] (4) to (\j);
      \foreach \j in {0,1}
      \draw [edge] (5) to (\j);
      \draw [edge] (0) to (1);
      \draw [edge] (2) to (3);
      \node (label) at (0,-1.3){$2K_3 + e$};
    \end{scope}
    
    \begin{scope}[xshift=0cm,yshift=-6cm]
      \foreach \i in {0,...,5}
      \node [vertex] (\i) at ({60*\i}:1){};
      \foreach \i in {0,...,4}
      \draw [edge] let \n1 = {int(mod(\i+1,6))}
      in
      (\i) to (\n1);
      \node (label) at (0,-1.3){$P_6$};
    \end{scope}
    
    \begin{scope}[xshift=2.5cm,yshift=-6cm]
      \foreach \i in {0,...,5}
      \node [vertex] (\i) at ({60*\i}:1){};
      \foreach \i in {0,...,3}
      \draw [edge] let \n1 = {int(mod(\i+1,6))}
      in
      (\i) to (\n1);
      \draw [edge] (3) to (5);
      \node (label) at (0,-1.3){$X_{172}$};
    \end{scope}
    
    \begin{scope}[xshift=5cm,yshift=-6cm]
      \foreach \i in {0,...,5}
      \node [vertex] (\i) at ({60*\i}:1){};
      \foreach \i in {0,...,3}
      \draw [edge] let \n1 = {int(mod(\i+1,6))}
      in
      (\i) to (\n1);
      \draw [edge] (3) to (5);
      \draw [edge] (4) to (5);
      \node (label) at (0,-1.3){$\overline{X_{98}}$};
    \end{scope}
    
    \begin{scope}[xshift=7.5cm,yshift=-6cm]
      \foreach \i in {0,...,5}
      \node [vertex] (\i) at ({60*\i}:1){};
      \foreach \i in {0,...,5}
      \draw [edge] let \n1 = {int(mod(\i+1,6))}
      in
      (\i) to (\n1);
      \node (label) at (0,-1.3){$C_6$};
    \end{scope}
  \end{tikzpicture}
  \caption{Graphs $J$ such that $J+P_3$ and $J+K_3$ are disconnected minimal
  \MT-obstructions.}
  \label{fig:obs-components}
\end{figure}
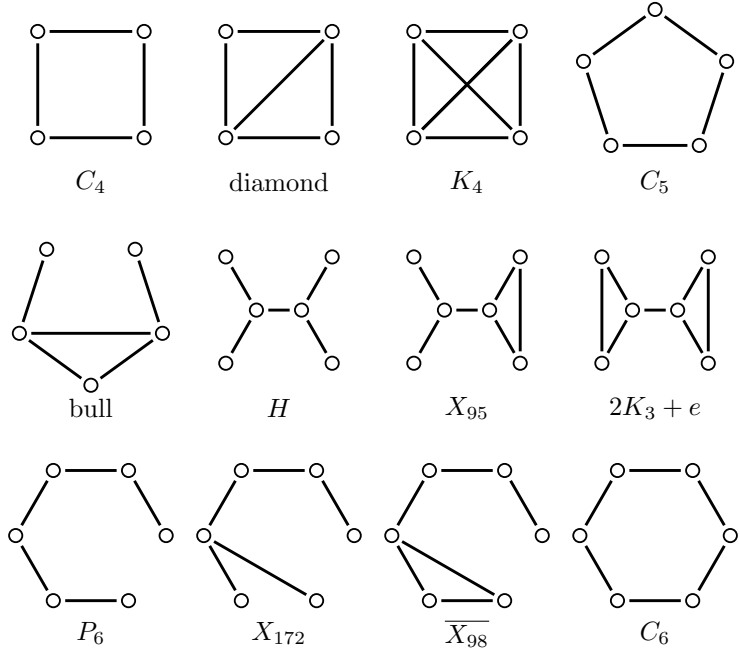

We start by studying the disconnected minimal \MT-obstructions.   

\begin{lemma}\label{lem:disconnected}
  Let $G$ be a disconnected \mayatupi graph. Then,
  \begin{enumerate}
    \item\label{item:no3components} $G$ has at most two components having at
      least three vertices; and
    \item\label{item:exactly2components} if $G$ has two components having $3$ or
      more vertices, then there must exist two non-adjacent vertices $u$ and $v$
      such that $G - \{ u,v \}$ is a disjoint union of copies of $K_1$ and
      $K_2$.  
  \end{enumerate} 
\end{lemma}
\begin{proof}
  Let $G$ be a disconnected \mayatupi graph and let $(A,B)$ be an
  $\MT$-partition of $G$. To prove \cref{item:no3components}, by contradiction
  suppose that $G$ has three connected components $C_1$, $C_2$ and $C_3$ having
  at least three vertices each. Since $\Delta(G[B])\leq 1$, note that at most 2
  vertices of $C_i$ lie in $B$, for every $i\in\{1, 2, 3\}$, and thus let
  $u_i\in V(C_i)\cap A$. Note that $\{u_1,u_2,u_3\}$ induce a $\overline{K_3}$
  in $G[A]$, a contradiction.
  
  Let us now prove \cref{item:exactly2components}. Assume that $G$ has exactly
  two components $C_1$ and $C_2$ having $3$ or more vertices and, similarly, let
  $u_i\in V(C_i)\cap A$, for each $i\in\{1,2\}$. Since $C_1$ and $C_2$ are
  connected components of $G$, note that no other vertex $w$ different from
  $u_1$ and $u_2$ can lie in $A$, as $w$ would be a neighbor of both $u_1$ and
  $u_2$, connecting $C_1$ to $C_2$. Consequently, $G-u_1-u_2 = G[B]$ is a graph
  formed by the union of copies of $K_1$ and $K_2$.
\end{proof}

With \cref{lem:disconnected} in mind, it is easy to observe that for each graph
$J$ in \cref{fig:obs-components}, the disjoint unions $J+P_3$ and $J+K_3$ result
in disconnected minimal \MT-obstructions with exactly two connected components.
By \cref{item:no3components} of \cref{lem:disconnected}, it is also clear that
$3P_3, 2P_3 + K_3, P_3 + 2K_3$, and $3K_3$ are disconnected minimal
\MT-obstructions.   Our following results show that these are the only
disconnected minimal \MT-obstructions.

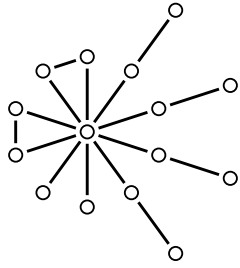
\begin{figure}[ht!]
  \centering
  \begin{tikzpicture}
    \node (c) [vertex] at (0,0){};
    \begin{scope}[rotate=-54]
      \foreach \i in {0,...,9}{
      \node (\i-1) [vertex] at ({36*\i}:1){};
      \draw [edge] (c) to (\i-1);
      }
      \foreach \i in {0,1,2,3}{
      \node (\i-2) [vertex] at ({36*\i}:2){};
      \draw [edge] (\i-1) to (\i-2);
      }
      \foreach \i/\j in {4/5,6/7}
      \draw [edge] (\i-1) to (\j-1);
    \end{scope}
  \end{tikzpicture}
  \caption{Nice non-complete components in the proof of \cref{lem:discon}.}
  \label{fig:discon}
\end{figure}

Define $\Fdisc$ to be the family of 28 disconnected graphs containing $3P_3,
2P_3 + K_3, P_3 + 2K_3$, $3K_3$ and all graphs obtained from the disjoint union
of a graph in \cref{fig:obs-components} with a $P_3$ or with a $K_3$. A graph
$G$ is \emph{nice} if it is obtained from a star by subdividing some edges once,
and then adding edges from some leaves to the central vertex (see
\cref{fig:discon}).

\begin{lemma}\label{lem:discon}  
  Let $G$ be an $\Fdisc$-free graph. If $G$ is disconnected and it has at least
  two connected components, $G_1$ and $G_2$, with at least three vertices each,
  then $G = n_1 K_1 + n_2 K_2 + G_1 + G_2$, for some $n_1,n_2$ non-negative
  integers. Moreover, $G_1$ and $G_2$ are nice.
\end{lemma}

\begin{proof}
  Note that $G$ cannot have three components with at least three vertices each
  from the fact that $3P_3, 2P_3 + K_3, P_3 + 2K_3, 3K_3 \in \Fdisc$.
  Consequently, $G = n_1 K_1 + n_2 K_2 + G_1 + G_2$, for some $n_1,n_2$
  non-negative integers.
  
  Let us prove that $G_1$ and $G_2$ are nice. If both $G_1$ and $G_2$ are
  complete graphs, then they are both isomorphic to $K_3$, because they have at
  least three vertices and $K_3 + K_4$ is in $\Fdisc$. Note that $G_1$ and $G_2$
  are nice.
  
  Hence, suppose that $G_2$ is not complete, so it contains an induced copy of
  $P_3$.   If $G_1$ is a complete graph, and since $K_4 + P_3$ is in $\Fdisc$,
  then it is a copy of $K_3$. In this case $G_1$ would be nice. The proof that
  $G_2$ is nice in this case is analogous to the one that follows assuming that
  $G_1$ is not complete as $J+K_3$ and $J+P_3$ are both in $\Fdisc$, for every
  $J$ in \cref{fig:obs-components}.
  
  Suppose that $G_1$ is not a complete graph. In this case, since $G_1$ and
  $G_2$ are not complete, proving that $G_2$ is nice is analogous to the proof
  that $G_1$ is nice.
  
  Since $C_4 + P_3, C_5 + P_3, C_6 + P_3$, and $P_6 + P_3$ are all in $\Fdisc$,
  we conclude that $G_1$ is a chordal graph.   Suppose first that $G_1$ is a
  tree.   If the diameter of $G_1$ is $2$, then $G_1$ is a star, and $G_1$ is
  nice.   If the diameter of $G_1$ is $3$, then $G_1$ is a \emph{double star},
  i.e. the graph obtained from two stars by adding an edge linking their
  centers. As $H + P_3$ is in $\Fdisc$ (see \cref{fig:obs-components} to verify
  the graph $H$), then one of the two centers of the double star has degree $2$,
  and thus, $G_1$ is nice. When the diameter of $G_1$ is $4$, there is a single
  \emph{center} $v$ in $G_1$, i.e. $v$ is the unique vertex whose distance to
  any other vertex in $G_1$ is at most two. Now, as $X_{172}$ is in $\Fdisc$,
  its absence in $G_1$ as an induced subgraph implies that the only vertex of
  degree possibly greater than $2$ in $G_1$ is $v$. Therefore, $G_1$ is as
  expected.
  
  We affirm that if $G_1$ is not a tree, then all triangles in $G_1$ share a
  common vertex (recall that $G_1$ is chordal and no longer cycles exist). It is
  impossible for $G_1$ to contain two triangles with disjoint sets of vertices,
  as otherwise $G$ would contain one of $2K_3 + P_3$, $(2K_3 + e) + P_3$, $C_4 +
  P_3$, $K_4 + P_3$, or $\textnormal{diamond} + P_3$ as an induced subgraph
  (depending on how many edges there are between them). Also, as
  $\textnormal{diamond} + P_3$ is in $\Fdisc$, two triangles in $G_1$ cannot
  share an edge.  Finally, the absence of $\textnormal{diamond} + P_3$ in $G$ as
  an induced subgraph prevents any three triangles in $G_1$ to pairwise share
  different vertices. Therefore, all triangles in $G_1$ share a vertex $v$.
  
  As $G$ is free from $C_4 + P_3$, $\textnormal{diamond} + P_3$, $K_4 + P_3$,
  and $\textnormal{bull} + P_3$, each triangle has a single vertex of degree
  greater than $2$, which is $v$.   But $G$ is also free from $X_{95}$ and
  $\overline{X_{98}}$, so the remaining vertices of $G_1$ are either pendant
  vertices attached to $v$, or belong to a pendant $P_3$ having $v$ as one
  extreme. Therefore, $G_1$ has the desired structure.
  
  As we observed before, recall that when $G_1$ is not a complete graph, we can
  apply the same arguments to conclude that $G_2$ also has the desired
  structure. When $G_1$ is a triangle and $G_2$ is not a complete graph, then
  all the above arguments still apply to verify that $G_2$ has the desired
  structure. This is because in each of the disconnected graphs in $\Fdisc$
  which were used in the proof, a $P_3$ component can be replaced by a $K_3$
  component, and the resulting graph is also in $\Fdisc$.
\end{proof}

\begin{theorem}\label{thm:disconnected}
  The disconnected graphs in $\Fdisc$ are exactly all the disconnected
  minimal \MT-obstructions.
\end{theorem}

\begin{proof}
  Clearly, disconnected minimal \MT-obstructions have neither isolated vertices
  nor connected components isomorphic to $K_2$. Thus, any disconnected
  \MT-minimal obstruction has at least two components with at least three
  vertices each.   But by virtue of \cref{lem:discon}, any such graph not in
  $\Fdisc$ admits an \MT-partition. Indeed, note that if $G = G_1+G_2$ and $G_1$
  and $G_2$ are nice, then one can place both centers of the stars that
  originated $G_1$ and $G_2$ in $A$, and the remaining vertices in $B$.
  Therefore, $\Fdisc$ contains every disconnected \MT-minimal obstruction.
\end{proof}

Recall that the class of \emph{cographs} is defined recursively as follows: a
trivial graph is a cograph; if $G_1, \dots G_k$ are cographs, then their
disjoint union $\sum_{i=1}^k G_i$ is a cograph; if $G$ is a cograph, then
$\overline{G}$ is a cograph~\cite{corneilDAM3}. They are also well-known to
correspond to the class of $P_4$-free graphs, and also to the class of graphs in
which every connected induced subgraph has a disconnected
complement~\cite{corneilDAM3}.

Being closed under taking complements, and being $P_4$-free, we can thus easily
determine minimal \MT-obstructions for cographs from \cref{thm:disconnected}.
Let $\Fcog$ be the family consisting of the graphs $3P_3$, $2P_3 + K_3$, $P_3 +
2K_3$, $3K_3, C_4 + J$, $\textnormal{diamond} + J$, and $K_4 + J$ for $J \in \{
P_3, K_3 \}$, and their complements.

\begin{corollary}
   A cograph $G$ is a \mayatupi graph if and only if it is $\Fcog$-free.
\end{corollary}

\begin{proof}
  Clearly, $\mathcal{F}_c$ is obtained by taking the intersection of all the
  disconnected minimal \MT-obstructions with the class of cographs, and then
  adding the complements of all the obtained graphs.   Since every connected
  cograph minimal \MT-obstruction has a disconnected complement, we rely on
  \cref{thm:disconnected} to conclude that $\mathcal{F}_c$ contains all
  cograph minimal \MT-obstructions.
\end{proof}


\section{Graphs with no induced four-vertex cycle}
\label{sec:c4-free}

In this section we present an algorithm to recognize \MT-graphs when restricted
to the class of $C_4$-free graphs in time $O(n^3)$.  An important ingredient of
our algorithm is a subroutine to efficiently recognize $C_4$-free $\langle 1,2
\rangle$-polar graphs; we begin this section by introducing such a subroutine.

In \cite{gagarinVNANBSFMN3}, Gagarin and Metelsky exhibited the complete list of
$\langle 1,2 \rangle$-polar minimal obstructions\footnote{This list can be
consulted in \url{https://graphclasses.org/classes/gc_664.html}}.   From this
list it is easy to conclude that the only chordal $\langle 1,2 \rangle$-polar
minimal obstructions are $2P_3$, $P_3 + K_3$, and $2K_3$.   Moreover, since
$C_5, C_6$ and $C_7$ are also $\langle 1,2 \rangle$-polar minimal obstructions,
it is clear that every $C_4$-free $\langle 1,2 \rangle$-polar graph is chordal.
Thus, we focus on the case of chordal graphs and use the list of chordal
$\langle 1,2 \rangle$-polar minimal obstructions to present a certifying
algorithm to verify whether a graph $G$ with vertex set $V$ and edge set $E$ is
a chordal $\langle 1,2 \rangle$-polar graph.   The proposed algorithm runs in
time $O(|V|+|E|)$.

\begin{algorithm}[!htbp]
\footnotesize{
    \SetAlgorithmName{Algorithm}{}
    \DontPrintSemicolon
    \SetKwData{False}{false}
    \SetKwData{True}{true}
    \SetKwFunction{New}{new}
    \SetKwFunction{Used}{used}
    \SetKwInOut{Input}{input}
    \SetKwInOut{Output}{output}
    \KwIn{A chordal graph $G$ with vertex set and adjacency lists ordered by
      degree non-increasingly.}
    \KwOut{A $\langle 1,2 \rangle$-polar partition $(A,B)$ of $G$ or a set of
      vertices inducing $2K_3, P_3 + K_3$ or $2K_3$}
    $K \gets$ a maximal clique constructed by iteratively choosing a candidate
      vertex of maximum degree\;
    $\{C_1, \dots, C_k\} \gets$ the connected components of $G - K$\;
    \If{$|V(C_i)| \le 2$ for each $i \in \{1,\dots,k\}$}{
        \Return $(K,\bigcup_{i=1}^k V(C_i))$\; \label{lin:all-good}
        }
    \If{$|V(C_i)|, |V(C_j)| \ge 3$ for $i,j \in \{1,\dots,k\}$}{
        $S \gets$ a $3$ vertex subset of $V(C_i)$ inducing a $P_3$ or a $K_3$\;
        $T \gets$ a $3$ vertex subset of $V(C_j)$ inducing a $P_3$ or a $K_3$\;
        \Return $S \cup T$\; \label{lin:all-bad}
        }
    {$C \gets$ the vertex set of the only component of $G-K$ with more than two
      vertices\;}
    \If{there is a vertex $u \in C$ with $|K|-1$ neighbours in $K$}{
        $x \gets$ the only non-neighbour of $u$ in $K$\;
        $S \gets$ a $3$ vertex subset of $C$ inducing a $P_3$ or a $K_3$\;
        \If{$x$ is adjacent to at least two components of $G-K$}{
            $T \gets x$ and two neighbours of $x$ in different components of
            $G-K$\;
            \Return $S \cup T$\; \label{lin:x-mid}
        }
        \If{$x$ is adjacent to a component $C_i$ of order at least $2$}{
            $T \gets x$ and $2$ vertices in $V(C_i)$ such that $T$ induces
            a $P_3$ or a $K_3$\;
            \Return $S \cup T$\; \label{lin:x-ext}
        }
        \If{there is a component $D$ of $G[C] - u$ with more than two
        vertices}{
            $S \gets$ a $3$ vertex subset of $V(D)$ inducing a $P_3$ or a
            $K_3$\;
            $y \gets$ a non-neighbour in $K$, different from $x$, of every
            vertex in $S$\;
            $w \gets$ a neighbour of $x$ not in $K$\;
            \Return $S \cup \{w,x,y\}$\; \label{lin:not-change-one}
        }
        $K \gets (K \setminus \{x\}) \cup \{u\}$\;
        \Return $(K, V(G) \setminus K)$\; \label{lin:change-one}
    }
    \If{there is a vertex $u \in C$ with $|K|-2$ neighbours in $K$}{
        $\{x,y\} \gets$ the two non-neighbours of $u$ in $K$\;
        \If{$x$ or $y$ has a neighbour $z$ outside $K$}{
            $S \gets$ a $3$ vertex subset of $C$ inducing a $P_3$ or a $K_3$\;
            \Return $S \cup \{x,y,z\}$\; \label{lin:miss-two-one}
        }
        \If{there is a component $D$ of $G[C] - u$ with more than two vertices}{
            $S \gets$ a $3$ vertex subset of $V(D) \setminus \{u\}$ inducing a
            $P_3$ or a $K_3$\;
            $z \gets$ a non-neighbour in $K$, different from $x$ and $y$, of
            every vertex in $S$\;
            \Return $S \cup \{x,y,z\}$\; \label{lin:miss-two-two}
        }
        $K \gets (K \setminus \{x,y\}) \cup \{u\}$\;
        \Return $(K, V(G) \setminus K)$\; \label{lin:change-two}
    }
    $S \gets$ a $3$ vertex subset of $C$ inducing a $P_3$ or a $K_3$\;
    $T \gets$ three vertices in $K$ non-adjacent to every vertex in $C$\;
    \Return $S \cup T$\; \label{lin:miss-three}
    \caption{\textsc{$\langle 1,2 \rangle$-polar Recognition}}
    \label{alg:12-tc-chordal}
    \DecMargin{1em}
}
\end{algorithm}

We begin with a simple technical lemma regarding cliques in chordal graphs.

\begin{lemma}
\label{lem:clique-poset}
    Let $G$ be a chordal graph and let $C$ be a clique in $G$.
    \begin{enumerate}
        \item For every pair of adjacent vertices $x$ and $y$ in $V \setminus C$
            we have $N(x) \cap C \subseteq N(y) \cap C$ or $N(y) \cap C
            \subseteq N(x) \cap C$.

        \item \label{item:comp-vertex} Let $S$ be a subset of $V \setminus C$
            such that $G[S]$ is connected.   If $C \subseteq  \bigcup_{s \in S}
            N(s)$, then there is vertex in $S$ that is complete to $C$.
    \end{enumerate}
\end{lemma}

\begin{proof}
    For the first item, if $v \in (N(x) \cap C) \setminus (N(y) \cap C)$ and $w
    \in (N(y) \cap C) \setminus (N(x) \cap C)$, then clearly $\{ v, w, x, y \}$
    induces a $C_4$ in $G$, which is impossible.
    
    Let $x$ be a vertex in $S$ having $|N(x) \cap C|$ as large as possible.   If
    $N(x) \ne C$, then $C \setminus N(x) \ne \varnothing$. Let $y$ be vertex in
    $S$ having a neighbour $v$ in $C$ which is not a neighbour of $x$, and
    choose $y$ closest to $x$ in $G[S]$ with this property.   Consider an
    $xy$-path $P$ of minimum length in $G[S]$.   The choice of $x$ implies that
    it has a neighbour $u$ in $C$ which is not a neighbour of $y$.   Suppose
    that $P = (x=w_0, \dots, w_k=y)$.  Since $y$ does not have $u$ as a
    neighbour and $w_{k-1}$ does not have $v$ as a neighbour, it follows from
    the first item that $(N(w_{k-1}) \cap C) \subseteq (N(y) \cap C)$ and $u$ is
    not a neighbour of $w_{k-1}$.   If $w_i$ is the last vertex of $P$ having
    $u$ as a neighbour, then $(u, w_i, \dots, y, v, u)$ is an induced cycle of
    length at least $5$ in $G$, contradicting the chordality of $G$.
    Therefore, $N(x) = C$.
\end{proof}

Although our following proposition is not necessary for the purposes of this
work, we think that it is of independent interest and the proof is short enough
to be included it here.   To the best of our knowledge, the equivalence of
chordality mentioned in \Cref{pro:chor-char} has not been published before.

Before, we need to introduce some terminology.   Given a graph $G$ and a subset
$S$ of its vertices is \emph{convex} in the \emph{geodetic convexity} if for
every pair of vertices $x,y$ in $S$, every vertex in every shortest $xy$-path in
$G$ also belongs to $S$. For further information on geodetic convexity, we refer
the reader to \cite{araujo2025introduction,pelayo2013geodesic}.

\begin{proposition}
\label{pro:chor-char}
    A graph $G$ is chordal if and only if for any clique $K$ of $G$ and any
    subset $S$ of $V \setminus K$ such that $G[S]$ is connected, if $K \subseteq
    \bigcup_{s \in S} N(s)$, then there is vertex in $S$ that is complete to
    $K$. Moreover, the set of all such vertices is a convex set of $G[S]$ in the
    geodetic convexity.
\end{proposition}
\begin{proof}
    For the necessity, \Cref{item:comp-vertex} of \Cref{lem:clique-poset} proves
    the existence of a vertex with the desired property.   Let $T$ be the subset
    of $S$ consisting of the vertices that are complete to $K$.   Suppose that
    there is a geodetic $P$ in $G[S]$ between two vertices $x$ and $y$ in $T$
    such that none of the internal vertices belongs to $T$.  Let $P'$ be the
    (possibly trivial) path $P - \{x,y\}$ in $G[S]$. Since $P'$ is connected and
    it is a subset of $V \setminus K$, the contrapositive of
    \Cref{item:comp-vertex} in \Cref{lem:clique-poset} implies the existence of
    a vertex $z$ in $K$ such that every vertex of $P'$ misses. Hence, $zxPyz$ is
    a cycle of length at least $4$ in $G$, contradicting the chordality of $G$.
    Therefore, $T$ is a convex set in $G[S]$.

    We prove sufficiency by contrapositive.   If $G$ contains an induced cycle
    $C$ of length greater than $3$, then by taking $K$ to be any edge of $C$,
    then the set $S$ consisting of the rest of the vertices of $G$ does not
    fulfill the required condition.
\end{proof}

The following lemma is the main ingredient of the anticipated algorithm to
recognize $C_4$-free $\langle 1,2 \rangle$-polar graphs.

\begin{lemma}
\label{lem:12-tc-alg-corr}
    \Cref{alg:12-tc-chordal} is correct.
\end{lemma}
\begin{proof}
    We verify that whenever \Cref{alg:12-tc-chordal} returns a set of vertices
    it induces either $2P_3, P_3+K_3$ or $2K_3$, and whenever it returns a
    partition of $V(G)$, it is a $\langle 1,2 \rangle$-polar partition.   We
    follow the notation of \Cref{alg:12-tc-chordal}.   Notice that since $K$ is
    a maximal clique, no vertex in $V(G) \setminus K$ is adjacent to every
    vertex of $K$, so the cases considered are exhaustive.
    
    The output is evidently correct when \Cref{alg:12-tc-chordal} finishes in
    \Cref{lin:all-good} or \Cref{lin:all-bad}.   It follows from
    \Cref{item:comp-vertex} in \Cref{lem:clique-poset} that when a vertex $u$ in
    $C$ with $|K|-1$ neighbours in $K$ exists, no vertex of $C$ is adjacent to
    the vertex $x$ that $u$ misses in $K$.   Hence, the set $T$ in
    \Cref{lin:x-mid,lin:x-ext} contains $x$ and vertices in components different
    from $C$, so $S \cup T$ induces one of $2P_3, P_3+K_3$ or $2K_3$.
    
    Suppose that \ref{alg:12-tc-chordal} finishes in \Cref{lin:not-change-one}.
    By the previous cases, we know that $x$ has at most one neighbour in $G-K$
    which lies in a trivial component of $G-K$.  So, the degree of $x$ is at
    most $|K|$.  When $K$ was constructed $x$ was considered before $u$ to be
    added to it, as otherwise $u$ would have been added to $K$ instead.   
    From here, we conclude that the degree of both $x$ and $u$ is $|K|$.   
    We claim that no other vertex in $C$ has $|K|-1$ neighbours in $K$. Indeed,
    if there is $u'$ in $C$ different from $u$ with $|K|-1$ neighbours in $K$,
    $u'$ has exactly one neighbour in $C$. In this case, since $C$ has at least
    three vertices and $u$ also has only one neighbour in $C$, there is a
    non-trivial $uu'$-path in $C$. The internal vertices of such a path either
    see $|K|-1$ vertices in $K$, contradicting the choice of $x$ when $K$ was
    constructed, or they see less than $|K|-1$ vertices in $K$, contradicting
    the chordality of $G$, as one may find an induced cycle of length at least
    4.  Hence, if a component $D$ of $G[C]-u$ with more than two vertices
    exists, it follows from \Cref{lem:clique-poset} that a vertex $y$ exists in
    $K \setminus \{x\}$ such that no vertex in $D$ is adjacent to $y$.
    Therefore, $\{x,y\}$ is anticomplete to vertices of $D$, and since the
    degree of $x$ is $|K|$, a vertex $w$ exists in neither $K$ nor $V(C)$ which
    is adjacent to $x$, so the set returned in \Cref{lin:not-change-one} induces
    the desired graph.

    If none of the previous cases occur, then the components of $(K \setminus
    \{x\}) \cup \{u\}$ are the components of $G-K$ except for $C$ and the
    component where $x$ lies.   But every component of $G[C]-u$ has at most two
    vertices, and the component to which $x$ belongs has two vertices.   So, the
    partition proposed in \Cref{lin:change-one} is a $\langle 1,2 \rangle$-polar
    partition.

    For the set returned in \Cref{lin:miss-two-one}, the previous case and
    \Cref{lem:clique-poset} imply that $z$ is not a vertex of $C$, so clearly
    $G[S]$ is anticomplete to $G[\{x,y,z\}]$.   Thus, $S \cup \{x,y,z\}$ induces
    one of $2P_3, P_3+K_3$ or $2K_3$.   If the algorithm finishes in
    \Cref{lin:miss-two-two}, then $x$ and $y$ have degree $|K|-1$.   By the
    construction of $K$, and because $x$ and $y$ were considered to join $K$
    before $u$, we conclude that $u$ also has degree $|K|-1$ and has exactly one
    neighbour in $C$.   As in the analysis of the execution finishing in
    \Cref{lin:not-change-one}, vertex $u$ is the only of $D$ having $|K|-2$
    neighbours in $K$, as otherwise there would be a vertex of degree greater
    than $|K|-1$ in $D$, or an induced cycle of length greater than $3$ would
    exist.   As both of the aforementioned cases are impossible,
    \Cref{item:comp-vertex} of \Cref{lem:clique-poset} implies the existence of
    a vertex $z$ which is anticomplete to $S$.   Therefore, the set returned in
    \Cref{lin:miss-two-two} induces a $P_3+K_3$ or a $2K_3$.

    If the algorithm reaches \Cref{lin:change-two}, then the vertex $u$ of $C$
    with $|K|-2$ neighbours in $K$ is such that every component of $G[C]-u$ has
    at most two vertices and neither $x$ nor $y$ has neighbours outside $K$.
    Hence, every component of $G - ((K \setminus \{x,y\}) \cup \{u\})$ has at
    most two vertices, and the proposed partition is a $\langle 1,2
    \rangle$-polar partition of $G$.

    Finally, if the execution of \Cref{alg:12-tc-chordal} finishes in
    \Cref{lin:miss-three}, every vertex in $C$ misses at least three vertices in
    $K$, so the algorithm chooses three such vertices together with three
    vertices in $C$ inducing a connected subgraph of $G$, and correctly returns
    these six vertices.
\end{proof}

All that remains to be done is a running time analysis of
\Cref{alg:12-tc-chordal}.

\begin{theorem}
    \Cref{alg:12-tc-chordal} is a $\mathcal{O}(n+m)$-time certifying algorithm
    to recognize $\langle 1,2 \rangle$-polar chordal graphs.
\end{theorem}
\begin{proof}
    \Cref{lem:12-tc-alg-corr} proves that \Cref{alg:12-tc-chordal} is a
    certifying algorithm to recognize $\langle 1,2 \rangle$-polar graphs when
    the input graph is chordal.   It is well known that each of chordal graph
    recognition, sorting the vertex set of a graph by degree, and sorting the
    adjacency lists of a graph with respect to a given ordering can be done in
    time (and space) $\mathcal{O}(n+m)$.   Thus, by receiving an arbitrary graph
    as input, is possible to preprocess it in the desired time to obtain a
    suitable input for \Cref{alg:12-tc-chordal}.

    The construction of the clique $K$ is also achievable in time
    $\mathcal{O}(n+m)$, as the neighbourhood of each vertex is explored once to
    test if it is possible to add it to $K$.   Once $K$ is created, it is
    possible to run BFS, also in time $\mathcal{O}(n+m)$, on $G-K$ to obtain its
    connected components.   In the last step, without changing the time
    complexity, it is also possible to verify the order of each component of
    $G-K$, and finish in one of \Cref{lin:all-good} or \Cref{lin:all-bad} if
    necessary.

    By exploring the neighbourhood of each vertex in $C$ a vertex with the
    maximum number of neighbours in $K$ can be found in time $\mathcal{O}(n+m)$,
    and the same time is sufficient to find vertices $x, y$ and $z$ used in the
    \no-certificates returned in
    \Cref{lin:not-change-one,lin:miss-two-one,lin:miss-two-two,lin:x-mid,lin:x-ext}.
    Finally, as vertices $x$ and $u$ have been already found for the
    aforementioned cases, it is possible to prepare the \yes-certificates
    returned in \Cref{lin:change-one,lin:change-two}  in constant time.
\end{proof}

We now turn our attention to $C_4$-free \mayatupi graphs.   Before arriving to
the main result of this section, we introduce a simple technical lemma that
partially describes the structure of a possible \MT-partition of a $C_4$-free
graph.

\begin{lemma}
\label{lem:c4free-mayatupi}
  If $G$ is a $C_4$-free \mayatupi graph and $(K,\barM,S,M)$ is an \MT-partition
  of $G$, then:
  \begin{enumerate}
    \item\label{item:c4freebarMatmosttwo} if $\barM \neq \varnothing$, then
      $|\barM| = 2$;
    \item\label{item:c4freeuniquepartition} if $\barM \neq \varnothing$ and
      $K \neq \varnothing$, then $G$ has a unique \MT-partition $(A,B)$ in which
      $A$ is maximal and $A$ is formed by the vertices $v$ such that $\barM
      \subseteq N_G(v)$.
  \end{enumerate}
\end{lemma}
\begin{proof}
  To prove \Cref{item:c4freebarMatmosttwo}, it suffices to remember that $\barM$
  induces a complete subgraph minus a perfect matching, thus it has an even
  number of vertices and the existence of at least four vertices in $\barM$
  implies the existence of an induced $C_4$ in $G$.

  Suppose now that $\barM \neq \varnothing$ and $K\neq \varnothing$. By
  \Cref{item:c4freebarMatmosttwo}, let $\{u,v\}=\barM$ and $w\in K$. For any
  vertex $z\in V(G)$ such that $u,v\in N_G(z)$, if $zw \notin E(G)$, then
  $\{u,v,z,w\}$ induces a $C_4$.  Therefore, since $G$ is $C_4$-free, we have
  $K\subseteq N_G(z)$.  Consequently, if $A$ is maximal, then $z \in A$. This
  proves \Cref{item:c4freeuniquepartition}
\end{proof}

As we mentioned before, the best time we know so far for recognizing \mayatupi
graphs and constructing an \MT-partition when one exists runs in time $O(n^8)$.
The main result of this section improves this time within the class of
$C_4$-free graphs to $O(n^3)$.

\begin{theorem}
\label{thm:c4free-alg}    
    \mtrecognition can be solved in $\mathcal{O}(n^3)$-time if the input graph
    $G$ is a $C_4$-free graph on $n$ vertices.
\end{theorem}
\begin{proof}
  First we verify whether $G$ has an \MT-partition $(A,B)$ such that $A$ has at
  most two vertices. This can be done in $\mathcal{O}(n^3)$-time by trying all
  subsets $A$ with at most two vertices in $V(G)$ and checking whether $B =
  V(G)\setminus A$ induces a graph such that $\Delta(G[B])\leq 1$.

  If no \MT-partition is found, then $G$ is \mayatupi if and only if it has an
  \MT-partition $(A,B)$ with $|A|\geq 3$. Let us first look for an \MT-partition
  in which $(K,\barM,S,M)$ satisfies $\barM \neq \varnothing$.

  By \cref{lem:c4free-mayatupi}, we have in this case that $|\barM| = 2$. Since
  $|A|\geq 3$, then $K\neq \varnothing$. Consequently, by
  \cref{lem:c4free-mayatupi}, for any possible subset $S=\{u,v\}$ with two
  vertices of $V(G)$, we verify if there is an \MT-partition $(K,\barM,S,M)$ in
  which $\barM = \{u,v\}$. It suffices to check whether the set $K$ of vertices
  $w \in V(G)$ such that $u,v\in N_G(w)$ induces a clique, and to check whether
  $V(G)\setminus K$ induces a subgraph with maximum degree at most one. This can
  be achieved in $\mathcal{O}(n^3)$-time.

  If no \MT-partition is found, by \cref{lem:c4free-mayatupi},
  \cref{item:c4freebarMatmosttwo}, the only possibility for $G$ to be \mayatupi
  is that $G$ has an \MT-partition $(A,B)$ in which $A$ induces clique. In this
  case, $G$ is \mayatupi if and only if $G$ is a $\langle 1,2 \rangle$-polar
  graph. This can be verified in brute-force thorough an $\mathcal{O}(n^6)$-time
  algorithm as Gagarin and Metelsky~\cite{gagarinVNANBSFMN3} proved that there
  are 18 minimal obstructions to this class each one with at most 6 vertices,
  but here we use our linear-time \Cref{alg:12-tc-chordal} to recognize chordal
  $\langle 1,2 \rangle$-polar graphs provided in as $C_4$-free $\langle 1,2
  \rangle$-polar graphs are actually chordal~\cite{gagarinVNANBSFMN3}.
\end{proof}

As highlighted in the last lines of the previous proof, \Cref{alg:12-tc-chordal}
seems to be an important tool in obtaining a better running time for the
recognition of $C_4$-free \MT-graphs. In this context, we propose the following
natural problem.

\begin{problem}
\label{prob:12tc}
  What is the optimal running time to decide whether a given graph $G$ is a
  $\langle 1,2 \rangle$-polar graph?
\end{problem}

A computational search yielded 54 minimal chordal $\MT$-obstructions having at
most 9 vertices. \Cref{thm:c4free-alg} shows a more efficient way to
check whether a chordal graph is \mayatupi.
\begin{corollary}
  Deciding whether a given chordal graph $G$ is \mayatupi can be done in
  $\mathcal{O}(n^3)$-time.
\end{corollary}

Our algorithm to recognize chordal \mayatupi graphs does not use the nice
properties of chordal graphs (e.g., their perfect elimination ordering).   So,
the following natural problem is proposed.

\begin{problem}
\label{prob:chordal}
  Find an efficient algorithm to decide whether a given chordal graph $G$ is
  \mayatupi.   In particular, can it be done in $\mathcal{O}(n+m)$-time?
\end{problem}


\section{Optimization problems}
\label{sec:optimization}

In this section, we investigate some classical optimization problems on
\mayatupi graphs. The {\emph{chromatic number}} of a graph $G$, denoted by
$\chi(G)$, is the least integer $k$ such that $G$ admits an assignment $c \colon
V(G) \to \{1, \ldots, k\}$ where for each edge $uv$ of $G$, $c(u) \neq c(v)$.
Deciding if a graph $G$ has $\chi(G) \leq k$ is one of Karp's 21 NP-complete
problems \cite{karp1972}. However, this parameter can be efficiently computed on
{\emph{perfect}} graphs which are graphs that do not contain induced odd cycles
or their complements \cite{grotschel1984}.

A lower bound for $\chi(G)$ is the {\emph{clique number}} of $G$, denoted by
$\omega(G)$, which is the size of the maximum clique of $G$. Originally, perfect
graphs were defined as those graphs $G$ such that $\chi(G) = \omega(G)$ and such
that this property holds for any of its induced subgraphs. There are two other
parameters closely related to the clique number of $G$: its {\emph{independence
number}}, which is the size of a maximum independent set of $G$, and its
{\emph{clique cover number}}, which is the minimum number of cliques that cover
the vertex set of $G$. For every graph $G$, its clique number and chromatic
number are equal to the independence number and clique cover number,
respectively, of its complement $\overline{G}$.  As the class of perfect graphs
is self-complementary, it turns out that any perfect graph $G$ has also its
clique cover number equal to its independence number, and this property holds
for any of its induced subgraphs.

Although \mayatupi graphs are a superclass of split graphs, which are perfect,
they admit $C_5$'s as induced subgraphs, and the determination of the parameters
above cited remains an interesting question to be approached. Notice that split
graphs are also chordal and so the {\emph{minimum fill-in number}}, which is the
minimum number of edges to be added in order to make a graph chordal, is zero.
Again, since \mayatupi graphs are not chordal, the study of this parameter is
also interesting in this class. A notion that is related to the fill-in is the
treewidth of a graph. We also investigate the treewidth of \mayatupi graphs.

In what follows, we provide polynomial-time algorithms to compute these
parameters in \mayatupi graphs. Let $G$ be a \mayatupi graph. Recall that an
$\MT$-partition $(K,\overline{M},S,M)$, where $K$ is maximum, of $G$ can be
obtained in in $\mathcal{O}(n^8)$ time. Our algorithms will work on a given
$(K,\overline{M},S,M)$ partition of $G$.

\subsection{Clique number and independence number}

Recall that a $(k,l)$-graph admits a partition of its vertex set into $k$
independent sets and $l$ cliques~\cite{brandstadtDM152}. Alekseev and
Lozin~\cite{AL2003} proved that computing a maximum weighted independent set of
a $(k,l)$-graph can be done in $\mathcal{O}(n^{6l+2})$-time, when $k\leq 2$.
Since \mayatupi graphs form a proper subclass of $(2,2)$-graphs, this result
implies that the weighted independence number of a \mayatupi graph can be
computed in $\mathcal{O}(n^{14})$-time. However, we can do better as we show in
the sequel.   Since the class of \mayatupi graphs is self-complementary,
withouth loss of generality we will discuss cliques instead of independent sets.

Notice first that there are three types of maximal cliques in a \mayatupi graph:
a vertex in $S \cup M$ together with all its neighbours in $K$ and one of its
neighbours for each non-edge in $\overline{M}$; two adjacent vertices in $M$
together with all their common neighbours in $K$ and one of their common
neighbours of each non-edge in $\overline{M}$; all the vertices in $K$ and one
vertex of each non-edge in $\overline{M}$.   The number of maximal cliques of
the third kind is $O(2^n)$, which in turn implies that for each vertex in $S
\cup M$, the number of cliques of the first two types is also $O(2^n)$.   Since
the number of vertices in $S \cup M$ is $O(n)$, the number of maximal cliques of
a \mayatupi graph is $O(n2^n)$.   Fortunately, we do not need to check all the
possibilities, we only want the one with the largest weight, so in practice, we
only have to check a linear number of cliques.

\begin{theorem}
  \label{thm:alphaomegaalg}
  Given a graph $G$ with $n$ vertices and $m$ edges, a weight function $w: V(G)
  \to \mathbb{R}^+$, and an $\MT$-partition $(K,\overline{M},S,M)$ of $G$, a
  maximum weighted independent set and a maximum weighted clique of $G$ can be
  computed in $\mathcal{O}(n+m)$-time.
\end{theorem}
\begin{proof}
  We just provide the algorithm to compute a maximum weighted clique of $G$. The
  algorithm to compute a maximum independent set of $G$ is similar.

  The algorithm starts by using time $O(n+m)$ to greedily color the vertices in
  $K \cup \overline{M}$ with colors $\{1, \dots, k\}$. Afterwards, it uses
  counting sort to order the vertices in these sets with respect to this
  coloring.   It is well-known that by using $O(n+m)$ space we can order the
  adjacency lists of all the vertices in $S \cup M$ in time $O(n+m)$ with
  respect to the ordering we just described.

  To calculate the maximum weighted clique of $G$, first, for each vertex $u$ in
  $S \cup M$, the algorithm calculates the maximum weight of a maximal clique
  including $u$.  This can be done in time $O(d(u))$ by simply adding the
  weights of $u$, the neighbours of $u$ in $K$, and the neighbour of $u$ with
  the largest weight for each non-edge of $\overline{M}$.

  Second, for each edge $uv$ in $M$, the algorithm calculates the maximum weight
  of a maximal clique including both $u$ and $v$.   This is achieved by
  traversing the adjacency lists of $u$ and $v$ at the same time to find their
  intersection, again, choosing each common neighbour in $K$ and the common
  neighbour with the largest weight for each non-edge of $\overline{M}$.   This
  step can done in time $O(d(u)+d(v))$ using two pointers because the orderings
  of the neighbourhoods of $u$ and $v$ are compatible.

  Third, the maximum weight of a maximal clique in $K \cup \overline{M}$ is
  computed, by choosing all the vertices in $K$ and the vertex with the largest
  weight for each non-edge in $\overline{M}$.
  
  Finally, the algorithm returns the maximum weighted clique set among all those
  computed.   Clearly, the third step uses $O(n)$ time, and the first two steps
  use time $\sum_{v \in V} O(d(v)) = O(m)$.   The desired time bound is then
  achieved.
\end{proof}

There are some differences to be taken into consideration for the maximum weight
independent set case, although the main strategy is very similar to the one we
used in \Cref{thm:alphaomegaalg}.   Start by calculating the weight of the
maximum independent set $I$ contained in $S \cup M$.   For each vertex $v$ in $K
\cup \overline{M}$, verify whether the weight of $v$ is greater than then sum of
the weights of its neighbours in $S \cup M$ (if a vertex sees both ends of an
edge in $M$ just consider the largest of the weights), if this test passes, then
the $N(v) \cap I$ substituted by $v$ in $I$.   Finally, the same can be done for
each non-edge $\{ u, v \}$ in $\overline{M}$; we need to calculate the
intersection of $N(u)$ and $N(v)$ to avoid double counting, but that can be done
efficiently as in the proof of \Cref{thm:alphaomegaalg}.   Clearly, a maximum
weight independent set must have been explored with this process.

\begin{corollary}
  \label{cor:alphaomegatotaltime}
  The weighted independence and weighted clique numbers of a vertex-weighted
  \mayatupi graph $G$ on $n$ vertices can be computed in
  $\mathcal{O}(n^8)$-time.
\end{corollary}
\begin{proof}
  Thanks to Proposition~\ref{pro:mt-recognition}, we can compute an
  $\MT$-partition of $G$ in time $\mathcal{O}(n^8)$. Given such a partition, the
  independence number and clique number of $G$ is then computed as in
  \Cref{thm:alphaomegaalg}.
\end{proof}

\subsection{Chromatic number and clique cover number}

In this section, we provide a polynomial-time algorithm to compute the chromatic
number and clique cover number of a \mayatupi graph. Since the chromatic number
of $G$ is equal to the clique cover number of $\overline{G}$, it suffices to
provide an algorithm to compute the chromatic number of $G$.

Our algorithm relies on a given $\MT$-partition $(K,\overline{M},S,M)$ of $G$.
We first propose a simple brute-force algorithm with worse time complexity, and
then we explain how to improve it. In general, as in the previous case of the
independence and clique numbers, the total running time of the algorithm is
$\mathcal{O}(n^8)$, $n$ being the number of vertices in the input graph $G$, as
we need to compute an $\MT$-partition of $G$ first.

\begin{proposition}
  \label{pro:chi-bounds}
  Given an $\MT$-graph $G$ with $n$ vertices and $m$ edges, and an
  $\MT$-partition $(K,\overline{M},S,M)$ of $G$, then $|K|+|\barM|/2\leq
  \chi(G)\leq |K|+|\barM|/2+2$.
\end{proposition}
\begin{proof}
  For the lower bound, if one removes exactly one vertex from each non-edge of
    $\barM$, then the remaining vertices of $\barM$ together with $K$ form a
    clique of size $|K|+|\barM|/2$. Hence, $\chi(G)\geq |K|+|\barM|/2$. For the
    upper bound, it suffices to color $K\cup \barM$ with $|K|+|\barM|/2$ colors
    (non-adjacent vertices in $\barM$ must receive the same color), and then to
    color $S\cup M$ with at most two additional colors. 
\end{proof}

\begin{theorem}
  \label{thm:chi-bruteforce}
  Given an $MT$-graph $G$ with $n$ vertices and $m$ edges, and an
  $\MT$-partition $(K,\overline{M},S,M)$ of $G$, an optimal proper coloring of
  $G$ can be computed in $\mathcal{O}(n^2m)$-time.
\end{theorem}
\begin{proof}
  Let $|K|+|\barM|/2 = k$. By \cref{pro:chi-bounds}, since $\chi(G)\leq k+2$,
  note that in any optimal coloring of $G$ at most two non-edges of $\barM$ can
  be colored in a way that their endpoints receive different colors. Thus, we
  can try all subsets $R$ of at most two vertices of $\barM$ to be those that
  will receive distinct colors from their non-neighbors in $\barM$. Then, we
  color $K\cup \barM$ accordingly, in linear time. By \textit{accordingly}, we
  mean that if $R=\varnothing$, then the vertices of $K$ receive distinct colors
  in $\{1,\ldots,|K|\}$, and the vertices of $\barM$ receive colors in
  $\{|K|+1,\ldots, k\}$, in such a way that they receive the same color as their
  non-neighbor in $\barM$. If $R = \{v\}$, the coloring is the same as before,
  except that $v$ is colored $k+1$, and if $R = \{v,w\}$, the coloring is the
  same as before, but with $v$ colored $k+1$ and $w$ colored $k+2$.
  
  Finally, we color $S\cup M$ greedily, using the least color in
  $\{1,\ldots,k+2\}$ that does not occur in its neighborhood in $K\cup \barM$.
  
  Concerning the correctness of the algorithm, since we check all possible
  subsets of $\barM$ with at most two vertices to be those that receive distinct
  colors from their non-neighbors, we are sure to consider all possible
  restrictions of an optimal coloring of $G$ to $K\cup \barM$. Consequently,
  since $S\cup M$ is formed by a disjoint union of $K_1$'s and $K_2$'s, being
  $K\cup \barM$ a separator of such components, we can color $S\cup M$ with at
  most two additional colors without creating any conflict with the coloring of
  $K\cup \barM$ through a greedy algorithm, and thus we are sure to find an
  optimal coloring of $G$.
  
  Since there are at most $\mathcal{O}(|\barM|^2)$ subsets of $\barM$ with at
  most two vertices, the time complexity of this algorithm is
  $\mathcal{O}(n^2(n+m))$.
\end{proof}

Let us now improve the previous algorithm to compute an optimal coloring of $G$
in $\mathcal{O}(nm)$-time. 

\begin{theorem}
  \label{thm:chi-improved}
  Given an $MT$-graph $G$ with $n$ vertices and $m$ edges, and an
  $\MT$-partition $(K,\overline{M},S,M)$ of $G$, an optimal proper coloring of
  $G$ can be computed in $\mathcal{O}(nm)$-time.
\end{theorem}
\begin{proof}
  First, without loss of generality, we can assume that $G$ is connected, as
  otherwise we can color each connected component separately. Consequently, note
  that $K\cup \barM\neq \varnothing$. Moreover, we assume that $K$ is
  \emph{maximal}, that is, each vertex in $S\cup M$ has a non-neighbor in $K\cup
  \barM$, as otherwise we could move it to $K$.

  Again, let $|K|+|\barM|/2 = k$. By \cref{pro:chi-bounds}, recall that $k\leq
  \chi(G)\leq k+2$. We first color $K\cup \barM$ with $k$ colors $\{1,\ldots,
  k\}$ as in \cref{pro:chi-bounds}. We shall consider a case analysis on the
  neighborhood of the vertices of $S\cup M$ in $K\cup \barM$ to find an optimal
  coloring of $G$ in each case. In some cases, we will need to recolor some
  vertices of $K\cup \barM$.

  Note that this preprocessing can be done in $\mathcal{O}(n+m)$-time, as we
  just need to detect the connected components of $G$, check the neighborhood of
  each vertex in $S\cup M$ and possibly move vertices to $K$, and define a
  partial coloring of $G$ as above restricted to $K\cup \barM$.

  We define that a non-edge $pq$ in $\barM$ \emph{distinguishes} two adjacent
  vertices of $uv$ in $M$ if one of $p$ and $q$ sees $u$ and misses $v$, and the
  other sees $v$ and misses $u$.

  \medskip\noindent{Case 1:} 
    There is a non-empty subset of edges $M_c$ of $M$ where each vertex contains
    a clique of size $k$ in its neighborhood in $K \cup \overline{M}$.

    In this case, note that each vertex in $M_c$ has all the colors in
    $\{1,\ldots, k\}$ in its neighborhood, but adjacent vertices of $M_c$ can be
    adjacent to distinct vertices of $\overline{M}$ having the same color. Also
    note that in this case $G$ has a clique of size $k+1$. Thus, $\chi(G)\in
    \{k+1,k+2\}$. In particular, $K$ is a subset of the neighborhood of all
    vertices in $M_c$. Note also that $\barM \neq \varnothing$, as otherwise we
    would contradict the maximality of $K$. Finally, notice that $M_c$ can be
    built in linear time by checking the neighborhood of each vertex of $M$ in
    $K\cup \barM$.
  
    We need to consider such adjacencies to find an optimal coloring of $G$. We
    distinguish two subcases depending on whether there is a non-edge of
    $\overline{M}$ that distinguishes all the adjacent vertices of $M_c$ or not.
    A non-edge $pq$ in $\barM$ is \emph{special} if it distinguishes all the
    adjacent vertices of $M_c$. 

    \medskip\noindent{Case 1.1:} For some edge $uv$ in $M_c$, $u$ and $v$ share
    a clique of size $k$ in their neighborhood in $K \cup \overline{M}$.

     In this case, $G$ has a clique of size $k+2$. Thus, we can color $G$ with
     $k+2$ as in the proof of \cref{pro:chi-bounds}, and we are done.

     Note that checking whether there is an edge $uv$ in $M_c$ such that $u$ and
     $v$ share a clique of size $k$ in their neighborhood in $K \cup
     \overline{M}$ can be done in $\mathcal{O}(n+m)$-time because $uv$ must be
     adjacent to all the vertices in $K$ and have a common neighbor for each
     non-edge of $\barM$. Completing the partial coloring of $K\cup \barM$ to an
     optimal coloring of $G$ can also be done in $\mathcal{O}(n+m)$-time, as in
     \cref{pro:chi-bounds}.

    \medskip\noindent{Case 1.2:} For each edge $uv$ in $M_c$, $u$ and $v$ do not
    share a clique of size $k$ in their neighborhood in $K \cup \overline{M}$.

     In this case, since each vertex in $M_c$ is adjacent to all the vertices in
     $K$ and to at least one vertex of each non-edge of $\barM$, observe that,
     for each edge $uv$ in $M_c$, there is at least one non-edge of
     $\overline{M}$ that distinguishes $u$ and $v$. We need to consider whether
     there is a non-edge of $\overline{M}$ that distinguishes all the adjacent
     vertices of $M_c$ or not, as we explain below.
    
    \medskip\noindent{Case 1.2.1:} There is a special non-edge in
    $\overline{M}$.

    To check whether $G$ has a special non-edge can be done in linear time,
    however checking both cases below is the most expensive part of our
    algorithm as, for each special non-edge of $\barM$, we need to check whether
    we can recolor one of its endpoints with color $k+1$ without creating any
    conflict with the vertices in $S\cup (M\setminus M_c)$. This can be done in
    $\mathcal{O}(nm)$-time.

    \medskip\noindent{Case 1.2.1.1:} There is a special non-edge $pq$ in
    $\overline{M}$ such that if we recolor $p$ with color $k+1$, then all the
    vertices in $M\setminus M_c$ and all the vertices in $S$ have a missing
    color in $\{1,\ldots, k+1\}$ in their neighborhood in $K\cup \barM$.

    In this case, we can color $G$ with $k+1$ colors as follows. Let $i$ be the
    color of $q$. First we recolor $p$ with color $k+1$. For the edges in $M_c$,
    we color those endpoints that are adjacent to $p$ with color $i$, and the
    ones adjacent to $q$ with color $k+1$ (recall that $p$ and $q$ distinguish
    all adjacent vertices of $M_c$). Finally, we color the vertices in $S\cup
    (M\setminus M_c)$ with a missing color in $\{1,\ldots, k+1\}$. Since $G$ has
    a clique of size $k+1$, such coloring is optimal.
    
    \medskip\noindent{Case 1.2.1.2:} For each special non-edge $pq$ in
    $\overline{M}$, if we recolor $p$ with color $k+1$, then there is a vertex
    in $M\setminus M_c$ or in $S$ that does not have a missing color in
    $\{1,\ldots, k+1\}$ in its neighborhood in $K\cup \barM$.

    In this case, we claim that $\chi(G) = k+2$. Suppose, for a contradiction,
    that $\chi(G) = k+1$. Let $c$ be a $(k+1)$-coloring of $G$. By the
    definition of $k$, note that there is at most one non-edge of $\overline{M}$
    whose endpoints receive distinct colors in $c$.
    
    Since $M_c\neq \varnothing$, such non-edge should exist, as otherwise all
    the vertices in $M_c$ would have the same set of $k$ colors in their
    neighborhood in $K\cup \barM$, and thus they should receive the same color
    in $c$, which is a contradiction.
    
    Let $pq$ be such a non-edge. Without loss of generality, assume that $c(q) =
    k$ and $c(p)=k+1$. If $pq$ is not an special non-edge, let $uv$ be an edge
    of $M_c$ not distinguished by $pq$. Since both vertices $u$ and $v$ are
    adjacent to a clique of size $k$ in $K\cup \barM$ and $pq$ does not
    distinguish $u$ and $v$, one of these vertices, say $u$, is adjacent to both
    $p$ and $q$. Consequently, all the colors in $\{1,\ldots,k+1\}$ appear in
    the neighborhood of $u$, a contradiction.

    Consequently, $pq$ is a special non-edge. By hypothesis, there is a vertex
    $x$ in $M\setminus M_c$ or in $S$ that does not have a missing color in
    $\{1,\ldots, k+1\}$ in its neighborhood in $K\cup \barM$ when we recolor $p$
    with color $k+1$. Since $pq$ is a special non-edge, $x$ is adjacent to both
    $p$ and $q$, and thus all the colors in $\{1,\ldots,k+1\}$ appear in the
    neighborhood of $x$, a contradiction.

    Thus, we can color $G$ with $k+2$ colors as in the proof of
    \cref{pro:chi-bounds}.

    \medskip\noindent{Case 1.2.2:} There is no special non-edge in
    $\overline{M}$.

    Let us prove that $\chi(G) = k+2$ in this case. Suppose, for a
    contradiction, that $\chi(G) = k+1$. Let $uv$ be any edge of $M_c$ and $pq$
    be any non-edge of $\overline{M}$ that distinguishes $u$ and $v$. Since
    there is no special non-edge in $\overline{M}$, there is an edge $u'v'$ in
    $M_c$ such that $pq$ does not distinguish $u'$ and $v'$. Analogously, let
    $p'q'$ be any non-edge of $\overline{M}$ that distinguishes $u'$ and $v'$. 
    
    Let $c$ be a $(k+1)$-coloring of $G$. By the definition of $k$, note that
    there is at most one non-edge of $\overline{M}$ whose endpoints receive
    distinct colors in $c$. Without loss of generality, assume that $p$ and $q$
    receive the same color in $c$. Note that $u$ and $v$ are have neighbors with
    the same set of $k$ colors in $K\cup\barM$, but then they should have the
    same color in $c$, which is a contradiction.
    
    Consequently, in this case we again color $G$ with $k+2$ colors in linear
    time as in the proof of \cref{pro:chi-bounds}. 

    \medskip\noindent{Case 2:} $M_c$ is empty, that is, there is no edge in $M$
    such that both of its endpoints see a clique of size $k$ in their
    neighborhood in $K \cup \overline{M}$.

    \medskip\noindent{Case 2.1:} There is a vertex in $S\cup M$ that is adjacent
    to a clique of size $k$ in $K \cup \overline{M}$.
    
    In this case, note that $G$ has a clique of size $k+1$. We claim that
    $\chi(G) = k+1$ in this case. To do this, for each vertex of $M$ that is
    adjacent to a clique of size $k$ in $K \cup \overline{M}$ (recall that there
    is at most one per edge of $M$), we color it with color $k+1$, as well as
    all vertices in $S$. The remaining vertices of $M$ are colored by  a missing
    color in the set $\{1,\ldots,k\}$ that does not appear in its neighborhood
    in $K\cup \barM$.

    \medskip\noindent{Case 2.2:} No vertex in $S\cup M$ is adjacent to a clique
    of size $k$ in $K \cup \overline{M}$.

    In this case, note that each vertex in $S\cup M$ has a color in $\{1,\ldots,
    k\}$ that does not appear in its neighborhood in $K\cup \barM$. Thus, in
    both cases below, we can color the vertices in $S$ with such a missing
    color. Let us treat vertices in $M$.

    \medskip\noindent{Case 2.2.1:} For each edge $uv$ in $M$, $u$ and $v$ have
    distinct missing colors in $\{1,\ldots, k\}$.

    In this case, we can color $u$ and $v$ with their missing colors, and thus
    we can color $G$ with $k$ colors. Recall that $G$ has a clique of size $k$.

    \medskip\noindent{Case 2.2.2:} There is an edge $uv$ in $M$ such that $u$
    and $v$ have the same missing color in $\{1,\ldots, k\}$.

    Let us prove that $\chi(G) = k+1$ in this case. Suppose, for a
    contradiction, that $\chi(G) = k$. Let $c$ be a $k$-coloring of $G$. By the
    definition of $k$, note that the coloring of $K\cup \barM$ is unique, the
    pairs of non-adjacent vertices of $\barM$ necessarily having the same color.
    Consequently, $u$ and $v$ have the same set of $k$ colors in their
    neighborhood in $K\cup \barM$, and thus they should receive the same color
    in $c$, which is a contradiction.

    Thus, for each edge of $M$ we color one of its endpoints with its missing
    color in $\{1,\ldots, k\}$, and the other endpoint with color $k+1$. In this
    way, we can color $G$ with $k+1$ colors.
    
    One may verify that Case 2 can also be verified in linear time, while a
    coloring of $G$ can be obtained in linear time as well. Thus, the total time
    complexity of the algorithm is $\mathcal{O}(nm)$.
\end{proof}

\subsection{Fill-in and treewidth number}

The {\emph{minimum fill-in}} of a graph $G$, denoted by $mfi(G)$, is the minimum
number of edges to be added to $G$ to turn it into a chordal graph. Let
$G=(V,E)$ be a connected non-chordal \mayatupi graph and $(K,\overline{M},S,M)$
one of its $\MT$-partition, where $K$ is maximal.

Consider the following procedure to fill-in $G$. For each non-simplicial vertex
$v$ of $S$, taking in account the new edges already included, add as many new
edges to $\overline{M}$ as needed in order to make $v$ simplicial. Observe that
the order in which the vertices of $S$ are considered does not change  
the number of edges to be added. Let us call this graph $G_1$. At the end of
this step, all the vertices of $S$ in $G_1$ are simplicial. 

Now let us examine the vertices of $M$. Let $uv$ be an edge of $M$. If $u$ and
$v$ are not neighborhood-comparable, and $|N(u)| \leq  |N(v)|$, then, by adding
edges incident to $v$ and neighbors of $u$ in $N(u) \setminus N(v)$ yields $N(u)
\subseteq N(v)$. Now, for every pair of non-adjacent vertices $\{x,y\}$ of
$\overline{M}$ seen by $u$ (respectively $v$), add the edge $xy$. Let us call
$G_2$ the resulting graph of this step.

Finally, suppose that $G_2$ has $l \geq 2$ non-edges in $\overline{M}$.
Arbitrarily choose one of them to not fill and fill-in all the others, in order
to obtain a graph $G_3$. Observe that the described procedure takes linear time
to be executed.

\begin{theorem}
Let $G=(V,E)$ be a non-chordal \mayatupi graph. Then the general algorithm
described above optimally fill-in $G$. 
\end{theorem}
\begin{proof}
    We can suppose that $G$ is connected since the algorithm can be executed
    independently in every of its components. Let us make some observations. 

    First, suppose that there is an edge $uv$ in $M$ such that $u$ and $v$ have
    non-comparable neighborhoods in $K \cup \overline{M}$, and let $q$ and $p$
    be vertices of $K \cup \overline{M}$ such that $\{p,u,v,q\}$ induces a $C_4$
    or $P_4$.  In the first case, it is mandatory to make the neighborhoods of
    $u$ and $v$ comparable. In the second case, either $\{p,u,v,q\}$ is part of
    a $C_5$ or a $C_6$, or $K=\varnothing$, $\overline{M} = \{p,q\}$, there is
    no vertex in $S \cup (M\setminus \{u,v\})$ adjacent to both $u$ and $v$, and
    no other pair of adjacent vertices in $M$ cover $\{p,q\}$. The last case is
    impossible, as $G$ would be already chordal.  Hence, $\{p,u,v,q\}$ is part
    of a $C_5$ or a $C_6$, and then again the neighborhoods of $u$ and $v$ must
    be made comparable. As a last observation, at least 2 (respectively 3)
    chords have to be included to make a $C_5$ (respectively a $C_6$) chordal.
    Moreover, if $\{p,u,v,q\}$ is part of $C_6$, there is another edge $u'v'$ of
    $M$ such that $u'$ and $v'$ have also incomparable neighborhood with respect
    to $q$ and $p$.  

    Second, since any two non-edges of $\overline{M}$ induce a $C_4$, any
    triangulation of $G$ must let at most one non-edge in $\overline{M}$.

    Third, observe that a perfect elimination ordering of $G_3$ can be obtained
    by enumerating $V(G_3)$ in this order: the vertices of $S$, $M$,
    $\overline{M}$ starting with eventually its unique non-edge, and finally
    $K$. Thus, $G_3$ is chordal.

    Now, consider an optimal fill-in $G_3'$ of $G$. We are going to argue that
    $G_3$, the chordal graph obtained by our algorithm, has at most as many
    edges as $G_3'$, and by consequence it is optimal. Considering the edges
    included to obtain $G_3$, we say that an edge is of type 1, if it was
    included to turn comparable the neighborhoods of two adjacent vertices of
    $M$, and of type 2, if its extremities are vertices of $\overline{M}$.
    Observe that every edge included to obtain $G_3$ is of type 1 or 2. 

    As we argued in the first observation above, given the partition where $K$
    is maximal, since $G$ is not chordal, for any edge $uv$ of $M$ and in any
    fill-in of $G$,  $N(u) \subseteq N(v)$ or $N(v) \subseteq N(u)$. Therefore,
    if $u$ and $v$ have incomparable neighborhoods, in $G_3'$ (as in $G_3$), we
    must have $N(u) \subseteq N(v)$ or $N(v) \subseteq N(u)$. The only way to do
    so is adding edges from $u$ or $v$ to, respectively, $N(v)$ or $N(u)$. Let
    $L= N(u) \cap N(v)$. If $|N(u)| \leq  |N(v)|$, according to our algorithm,
    $|N(u)| - |L|$ edges will be added to make $N(u) \subseteq N(v)$, which is
    less than or equal to $|N(v)| - |L|$, in the case we turn $N(v) \subseteq
    N(u)$. By consequence, the number of edges of type 1 included to produce
    $G_3$ is minimum and $G_3'$ has at least this number of edges of type 1
    included to get the necessary neighborhood inclusion.

    Now, let us examine the inclusion of edges of type 2. Because of the second
    observation, regarding this type of edges, the only possibility of
    difference between $G_3'$ and $G_3$, in number of edges, is when all the
    edges between the pairs of $\overline{M}$ are included, instead of
    $|\overline{M}|-1$ edges, which have to be necesssarily included.  Suppose
    that $p,q$ is a pair of $\overline{M}$ which is not an edge in $G_3'$. If
    $pq$ was included to obtain $G_1$, then there is a  vertex $v \in S$
    adjacent to both $p$ and $q$. Since $K$ is maximal and $G$ is not chordal,
    $K \cup \overline{M} \neq \{p,q\}$ and $v$ misses a vertex $r$ of $K \cup
    \overline{M} \setminus \{p,q\}$. Therefore $\{v,p,q,r\}$ induces a $C_4$ in
    $G$. So, $G_3'$ contains the edge $vr$ and $G_3$ not, since $vr$ is not of
    type 1 or 2. Nevertheless, $G_3$ has at most as many edges as $G_3'$.

    Now, suppose that $pq$ was included in the process to obtain $G_2$. Then,
    there is an edge $uv$ of $M$, such that $u$ or $v$, let us suppose $v$, is
    adjacent to both $p,q$. If it is already the case in $G$, then by the same
    argument previously used, $v$ misses a vertex $r$ of $K \cup \overline{M}
    \setminus \{p,q\}$, $vr \in E(G_3')$ and not in $G_3$, since it is not of
    type 1 or 2, implying that $G_3$ has at most as many edges as $G_3'$. So,
    suppose that the edge $vq \notin E(G)$ and it was included to make
    comparable the neighborhoods of $u$ and $v$. Then $\{u,v,p,q\}$ induces a
    $P_4$ of $G$, and  as previously argued in the first observation, it makes
    part of $C_5$ or a $C_6$. Let $r$ be any vertex of $K \cup \overline{M} \cup
    M$ such that $\{u,v,p,q,r\}$ induces a $C_5$. Since $u$ and $v$ have
    comparable neighborhoods in $G_3'$, we can suppose that $uq \in E(G_3')$.
    Now, since $pq \notin E(G_3')$, and $\{u,p,r,q\}$ is a $C_4$, $ur \in
    E(G_3')$. Observe that $ur \notin E(G_3)$ since it is not of type 1 or 2.
    Again, we get that $G_3$ has at most as many edges as $G_3'$. On the other
    hand, if $\{u,v,p,q\}$ is a part of a $C_6$, then there is another edge
    $u'v'$ of $M$, such that $u'$ and $v'$ have incomparable neighborhoods and
    such that $\{u',v',p,q\}$  induces a $P_4$. Again, $G_3'$, as $G_3$,
    contains exactly one of $u'q$ or $v'p$, let us suppose $u'q$. Now, $\{u, q,
    u',p\}$ is a $C_4$ in $G_3'$ (and in $G_3$).  Since $pq \notin E(G_3')$,
    $uu' \in E(G_3')$ and not in $G_3'$ since it is not of type 1 or 2.
    Therefore, again $G_3$ has at most as many edges as $G_3'$.

    The last step of the algorithm, to obtain $G_3$ is a direct consequence of
    the second observation we have done before. Since at most one non edge of
    $\overline{M}$ can last unfilled, we have added a minimum number of edges.
    Likewise $G_3'$, $G_3$ is chordal and has at most the its same number of
    edges and so is an optimal fill-in of $G$.
\end{proof}

There is a close relationship between a fill-in and a tree decomposition of a
graph $G$.

A \textit{tree decomposition} of a graph $G$ is a pair $(T,\mathcal{B}=\{B_t
\mid t \in V(T)\})$, where $T$ is a tree and each set $B_t$, called a
\textit{bag}, is a subset of $V(G)$, satisfying the following properties: 
\begin{enumerate}
    \item $\bigcup_{t \in V(T)}B_t = V(G)$,
    \item for every edge $uv \in E(G)$, there exists a bag $B_t$ with $u,v \in
    B_t$, and
    \item for every vertex $v \in V(G)$, the set $\{t \in V(T) \mid v \in B_t\}$
    induces a connected subgraph of~$T$. 
\end{enumerate}
The \textit{width} of a tree decomposition $(T, \{B_t \mid t \in V(T)\})$ is
$\max_{t \in V(T)}|B_t|-1$, and the \textit{treewidth} of a graph $G$, denoted
by $\tw(G)$, is the minimum width of a tree decomposition of $G$.

Given a tree decomposition $(T, \{B_t \mid t \in V(T)\})$ of $G$, it is
well-known that by filling-in each bag of $T$, i.e. by adding all missing edges
among vertices within each bag, a chordal graph $G'$ is produced. Consequently,
each bag will correspond to a clique of such chordal supergraph $G'$ and the
treewidth of $G$ is the size of a largest clique of $G'$ minus one. Indeed, this
is an equivalent characterization of the treewidth of a graph $G$ as being the
size of largest clique in a triangulation of $G$ minus one.

Given a fill-in $G'$ of $G$, it is clear that $\omega(G') - 1$ is an upper bound
for the treewidth of $G$. However, it is not always the case that when we try to
minimize the number of edges to be added to $G$ to make it chordal, we also
minimize the size of the largest clique of the resulting chordal graph (see
\cite{ghahremani2022chordal,Dereniowski2019}).

We now prove that the triangulation of $G$ given by our algorithm produces a
graph with the minimum size of the largest clique as possible and therefore, the
treewidth of $G$ is obtained as a byproduct.

In the proof of the next theorem, we use the following well-known remarks:
\begin{lemma}\ 
  \label{lem:tw}
  \begin{enumerate}[(a)]
    \item For every graph $G$ and every subgraph $H$ of $G$, $\tw(H) \leq
    \tw(G)$.
    \item If $H = (A\cup B, E)$ is a complete bipartite graph and $H\subseteq
    G$, then for any tree decomposition $(T, \cal{B})$ of $G$, there exists a
    bag $B_t\in \mathcal{B}$ such that $A\subseteq B_t$ or $B\subseteq B_t$.
    \item If $K$ is a clique of $G$, then for any tree decomposition $(T,
    \cal{B})$ of $G$, there exists a bag $B_t\in \mathcal{B}$ such that
    $K\subseteq B_t$.
  \end{enumerate}
\end{lemma}

\begin{theorem}
    Let $G$ be a non-chordal \mayatupi graph and $G_3$ be the graph produced by
    the general algorithm described above. Then, $tw(G)=tw(G_3)=\omega(G_3) -
    1$.  
\end{theorem}
\begin{proof}
    Let $(K,\overline{M},S,M)$ be an $\MT$-partition of $G$ where $K$ is
    maximal. Let $(T, \mathcal{B} = \{B_t \mid t \in V(T)\})$ be a optimum tree
    decomposition of $G$ and $G'$ be the triangulation of $G$ given by this tree
    decomposition, i.e. $G'$ is obtained from $G$ by adding edges linking
    pairwise non-adjacent vertices that lie in some $B_t\in \mathcal{B}$. Let $L
    = \{e_1, e_2, e_3, \dots, e_l\}$ be a non-empty set of non-edges in
    $\overline{M}$. We claim that $T$ has a bag containing $K \cup (\overline{M}
    - e_i)$, for some edge $e_i$ of $L$. Let $e_i = \{p,q\}$ be any non-edge of
    $L$. Since $p$ and $q$ are non-adjacent and have $K \cup (\overline{M} -
    \{p,q\})$ as common neighbors, there is bag that contains  $K \cup
    (\overline{M} - \{p,q\})$ or a bag that contains $\{p,q\}$, by
    Lemma~\ref{lem:tw}. In the second case, the triangulation of $T$ will
    include the edge $pq$. Since the same occurs for every non-edge $e_i$ of
    $L$, then $T$ contains a bag $X$ with $K$ and all the pairs of $L$, except
    maybe one pair $e_i$. Therefore, $G'$ contains all, except maybe one, the
    edges $L = \{e_1, e_2, e_3, \dots e_l\}$. Suppose now that for every
    non-edge $e_i$ of $L$, its extremities $p$ and $q$ have a common neighbor $s
    \in S \cup M$. Our algorithm will turn $L$ into edges producing a clique $K
    \cup \overline{M}$ in $G_3$. Suppose this edge do not exist in $G'$, which
    means that $p$ and $q$ do not lie in a same bag of $T$. Since $\{s\} \cup K
    \cup (\overline{M} - \{p,q\})$ are common neighbors of $p$ and $q$, then
    there is a bag containing $\{s\} \cup K \cup (\overline{M} - \{p,q\})$ by
    Lemma~\ref{lem:tw}, which results in a clique of size $1+ |K \cup
    (\overline{M} - \{p\})|$ in $G'$, the same size of the clique produced by
    our algorithm, but in our case, by the inclusion eventually of less edges.
    Additionally, if $p$ and $q$ are seen by more than two vertices of $S \cup
    M$, then adding the edge $pq$ will produce a smaller clique than turning the
    neighborhood of $p$ and $q$ into a clique. By consequence, turning into
    edges all the pairs of vertices listed in $L$ having a common neighbor
    outside $K \cup \overline{M}$ does not produce bigger cliques than those
    existent in $G'$. Moreover, at most one pair of vertices of a non-edge in
    $\overline{M}$ cannot belong to the same bag. 

    Now, let us examine the edges included in order to produce neighborhood
    inclusion between two vertices $u$ and $v$ such that $uv$ is an edge of $M$.
    Let us examine $T$ again and its triangulation producing $G'$. Since $uv$ is
    an edge of $G$, the neighborhoods of $u$ and $v$ should be comparable,
    otherwise $G'$ would have a $C_4$. Hence, either $N_{G'}(u) \subseteq
    N_{G'}(v)$ or $N_{G'}(v) \subseteq N_{G'}(u)$. Observe that $N_{G'}(u) \cup
    N_{G'}(v) \subseteq X$, maybe except by a pair of vertices. Therefore,
    producing neighborhood inclusion will create an almost clique $N_{G'}(u)
    \cup N_{G'}(v) \cup \{v\}$ or $N_{G'}(u) \cup N_{G'}(v) \cup \{u\}$ in $G'$,
    where an almost clique is a set of vertices inducing a complete graph less
    one edge. They have clearly the same size, which is the size of the almost
    clique produced by our algorithm, but in our case, adding a minimum number
    of edges. If $N_{G'}(u) \cup N_{G'}(v)$ do not belong to $X$, as we have
    argued before, at most one non-edge of $\overline{M}$ seen by $u$ and $v$ do
    not lie in $X$. But again, including this pair in $X$ will produce a clique
    of size no greater that the clique containing the common neighbors of this
    non-edge. 

    As a conclusion, since $G_3$ has no larger cliques than $G'$, they have the
    same size of maximum clique, which is equal to $tw(G) + 1$. 
\end{proof}


\section{Conclusions and open problems}
\label{sec:conclusions}
In this work, we introduce a new, to the best of our knowledge, class of graphs
that properly contain the well-known split graphs. Since it falls in the
sparse-dense scheme~\cite{federSIDMA16}, we can recognize a graph in this class
in $\mathcal{O}(n^8)$ time. For the class of split graphs, a proper subclass, we
know how to recognize them in linear time. The main open problem we have is:

\begin{problem}
    \label{prob:linear-recog}
    Can one recognize a \mayatupi graph in linear time?
\end{problem}

We show several particular cases with a better algorithm is provided, such as
when $G$ has bounded clique or independence numbers, or when $G$ is $C_4$-free,
or when $G$ has bounded neighborhood diversity. For the classes of \mayatupi
graphs that are cographs, or trees, we not only exhibit improved recognition
algorithms, but we also provide a complete characterization by forbidden induced
subgraphs. Since we obtained with computer assistance over than 2000 minimal
obstructions for the class of \mayatupi graphs, we also ask:

\begin{problem}
  \label{prob:obstructions}
    Can one find a complete list of minimal obstructions for the class of
    chordal \mayatupi graphs?
\end{problem}

Finally, we have shown that the chromatic number of a \mayatupi graph can be
computed in polynomial time. We have also provided a polynomial time algorithm
to compute an optimal fill-in of a \mayatupi graph, and as a byproduct, to
compute its treewidth. 

These problems are known to be linear time solvable for split graphs, however,
to the best of our knowledge\footnote{This is also reported as unknown in the
\url{graphclasses.org} database.}, there are no efficient algorithms proposed to
solve these problems for (2,2)-split graphs, which are a proper superclass of
\mayatupi graphs. Therefore, we ask: 

\begin{problem}
  \label{prob:complexity}
  Can one find a decision problem that is in $\P$ for split graphs, but is
  $\NP$-complete for \mayatupi graphs?
\end{problem}

\section{Acknowledgments}

We are deeply thankful to Professor Andrei Gagarin for sharing a scanned copy of
\cite{gagarinVNANBSFMN3} with us.

J\'ulio Ara\'ujo is partly funded by CNPq (Brazil) projects
404479/2023-5, 404613/2023-3 and 308939/2025-5, Inria Associated Team CANOE
(Brazil-France) and project CAPES-Cofecub
49587PE (Brazil-France). César Hernández-Cruz is supported by UNAM-PAPIIT IN106425 and
SEP-CONACYT CB A1-S-8397. Cl\'audia Linhares Sales is supported by projects number
312044/2022-4 and 404479/2023-5 of CNPq - Conselho Nacional de Desenvolvimento
Científico e Tecnológico, as well as the Program PREI-DGAPA of the Universidad
Nacional Aut\'onoma de M\'exico.

\vspace{2mm}
\bibliographystyle{acm}
\bibliography{manuscript}

\end{document}